\documentclass[12pt,twoside,english,a4paper]{article}
\usepackage{etex}

\usepackage{caption}
\usepackage[latin1]{inputenc}
\usepackage[intlimits] {amsmath}   % muss for defjs# stehen
\usepackage{graphicx}
\usepackage{psfrag}
\usepackage{ifthen}
\usepackage{fancyhdr}
\usepackage{rotating}
\usepackage{multirow}
\usepackage{booktabs}              % Dicke Tabellenlinien
\usepackage{amssymb}
\usepackage{amsbsy}
\usepackage{bm}                    % fette Schrift in Formeln
\usepackage{bbm}
\usepackage{babel}
\usepackage{theorem}
\usepackage{upgreek}  % gerade (roman) griechische Buchstaben z.B. \upalpha
\usepackage{pifont}   % zusaetzliche Schriften, z.B. eingekreiste Zahlen
\usepackage[active]{srcltx}   % aktive suche im xdv

\usepackage{suetterl}  % Suetterlin Schrift
\newfont{\suetdbl}{suet14 scaled 2000}  % Suetterlin skaliert um Faktor 2,000

\usepackage{yfonts}  % altdeutsche Schrifte und Initialen
\newfont{\gothdbl}{ygoth scaled 2000} 
\newfont{\frakdbl}{yfrak scaled 2000}
\newfont{\swabdbl}{yswab scaled 2000}
\usepackage[normalem]{ulem}  % zusaetzliche Unterstreichungsformen wie z.B.:
\usepackage[textsize=footnotesize,textwidth=1.7cm]{todonotes}
\usepackage[]{defjs5}      % eigene Definitionen: Option BoldVarsRoman

\usepackage{float}
\usepackage{caption}
\usepackage{wrapfig}
\usepackage[authoryear]{natbib}
\bibpunct{[}{]}{,}{n}{}{;}
\fancypagestyle{plain}{%
  \fancyhead{}%
  \fancyfoot[c]{\sffamily\thepage}%
}
\makeatletter                           % Leere Fuellseiten
\def\cleardoublepage{\clearpage\if@twoside \ifodd\c@page\else
  \hbox{}
  \vspace*{\fill}
  \thispagestyle{empty}
  \newpage
  \if@twocolumn\hbox{}\newpage\fi\fi\fi}
\makeatother
\begin{document}
\unitlength1.0cm
\frenchspacing

%=== cover page
%\include{fContact_03}

%\tableofcontents
\newpage

%=== sections ===
%--- section 1
\thispagestyle{empty}
\ce{\bf \large
A robust mixed finite element formulation for
}
\ce{\bf \large
third medium contact
}

%\vspace{-3mm}
%\begin{center}
%{\bf \large
%Subtitle of the report.
%}
%\end{center}

\vspace{2mm}
\ce{M. Vorwerk$^1$, J. Schr\"oder$^1$ and P. Wriggers$^2$}

\vspace{2mm}
\ce{$^1$Institute of Mechanics, Faculty of Engineering, University of Duisburg-Essen}
\ce{\small e-mail: maximilian.vorwerk@uni-due.de, j.schroeder@uni-due.de}
\vspace{2mm}
\ce{$^2$Institute of Continuum Mechanics, Faculty of Mechanical Engineering, Leibniz University Hannover}
\ce{\small e-mail: wriggers@ikm.uni-hannover.de}

\vspace{2mm}
\begin{center}
{\bf \large Abstract}
\bigskip

{\footnotesize
\begin{minipage}{16.cm}
\noindent
Third medium contact provides a smooth continuum alternative to classical contact algorithms by replacing explicit contact constraints with a highly compliant fictitious medium.  
In this work, an auxiliary-field stabilization is introduced in which a deformation-gradient-like field is treated as an independent unknown in the third medium and coupled to the physical deformation gradient by a penalty term.  
A gradient contribution acting on the auxiliary field provides the regularization mechanism without requiring a direct evaluation of higher displacement derivatives.  
Linear and quadratic interpolation spaces are investigated, including continuous and element-wise discontinuous auxiliary-field approximations.  
The numerical results show that continuous low-order auxiliary fields provide an effective gradient-type stabilization of the third medium, even when the displacement field is approximated by first-order finite elements.  
For element-wise discontinuous auxiliary fields, the additional unknowns remain local to each element and can be eliminated locally by static condensation, so that the global system does not necessarily contain additional auxiliary degrees of freedom.  
Benchmark problems involving large deformation, progressive self-contact and severe third-medium compression are used to assess the formulation.

\end{minipage}
}
\end{center}

%{\bf Keywords:}
%third medium contact, finite deformation, auxiliary-field stabilization, mixed finite elements, low-order finite elements

\vspace{-5mm}
%=========================================================================
\section{Introduction}
\vspace{-5mm}
%=========================================================================
%
The finite element method is a central computational tool of virtual prototyping and simulation-driven engineering design.
It enables simulation-driven design before physical prototypes are built and motivates robust numerical methods for strongly nonlinear boundary value problems.
Contact mechanics is a prominent example. 
Its difficulty arises from evolving contact zones, changing boundary conditions, large local deformations and localized force transfer, cf.~\cite{Wri:2006:ccm,Lau:2002:cca,Yas:2013:nmi}. 
These effects occur in metal forming, crash simulations and sealing, but also in soft robotics and contact-aided mechanisms, where contact becomes part of the desired response, cf.~\cite{RusTol:2015:dfa,FreSigPou:2024:too}.
Classical finite element contact formulations discretize potential contact interfaces explicitly and enforce the constraints by penalty, Lagrange multiplier, augmented Lagrangian or barrier-type methods, cf.~\cite{Wri:2006:ccm,Lau:2002:cca}. 
Surface-to-surface and mortar formulations extend this framework to non-matching contact meshes, cf.~\cite{PusLau:2004:ams,SauDeL:2013:acc}.
Despite their maturity, these methods rely on contact surfaces, gap functions, active sets and contact search algorithms. 
For large sliding, self-contact, severe mesh distortion or topology optimization, where contact boundaries may evolve or may not be known in advance, this becomes a major algorithmic burden, cf.~\cite{BluSigPou:2021:icm,FreSigPou:2024:too,WriKorJun:2025:atm}. 
This motivates formulations in which contact emerges from a continuum model instead of being imposed on an explicit lower-dimensional interface.
The third medium contact method follows this idea. 
Instead of enforcing contact constraints directly between two surfaces, the space between potentially contacting bodies is filled with a highly compliant fictitious medium. 
This third medium has only a negligible influence before contact, but stiffens strongly when compressed to nearly vanishing volume. 
Contact forces are then transmitted through the compressed medium, while explicit contact search and inequality constraints are avoided. 
The method is introduced in \cite{WriSchSch:2013:afe} for frictionless finite deformation contact using an isotropic-anisotropic material model for the intermediate medium. 
A related fictitious contact material combined with high-order finite elements is investigated in \cite{BogZanKolRan:2015:ncw}. 
The extension to isogeometric analysis, including a study of material parameters and higher-order spatial convergence, is presented in \cite{KruNguWriDeL:2018:ifc}. 
An isogeometric-meshfree coupling strategy is proposed in \cite{HuaNguZho:2018:aim}, showing that the concept can also be transferred beyond standard finite element discretizations.
A renewed interest in third medium contact emerges in density-based topology optimization. 
In this setting, the void phase naturally provides a region in which a fictitious contact medium can be placed, while contact boundaries are not known in advance. 
Internal contact modeling for finite strain topology optimization is introduced in \cite{BluSigPou:2021:icm}. 
This work also proposes the so-called HuHu regularization, which penalizes the Hessian of the displacement field in the third medium to control severe element distortions. 
The use of third medium contact for topology optimization of self-contacting structures is further developed in \cite{FreSigPou:2024:too}, where additional design requirements are introduced to improve the robustness of optimized contact-aided mechanisms. 
Internal contact is also used as a design mechanism for nonlinear elastic springs with tailored force-displacement responses in \cite{BluSigPou:2023:ido}. 
The concept is extended to contact-aided thermo-mechanical regulators in \cite{DalAleFrePouSig:2025:too}, where self-contact is used to tune heat transfer through switches, diodes and triodes. 
Frictional third medium contact is introduced in \cite{FreRokPouSigGee:2024:aft} by adding an anisotropic shear contribution and a crystal-plasticity-inspired slip mechanism that mimics Coulomb friction.

Despite these advances, the stabilization of the third medium remains a central issue. 
The medium must be soft enough to avoid forces before contact, but stable enough to withstand extreme compression, shear and mesh distortion during contact. 
The HuHu regularization is highly effective, but it penalizes all second-order deformation modes, including bending. 
To reduce this effect, the HuHu-LuLu regularization is introduced in \cite{FreDalSigPou:2025:itm}. 
There, a Laplacian contribution is subtracted from the Hessian contraction, which reduces the penalization of bending and quadratic compression while maintaining control over undesirable skew deformation modes. 
A deformation-gradient averaging regularization is proposed in \cite{FalAmaHor:2026:dga}. 
This approach penalizes deviations of the deformation gradient at the integration points from an element-wise representative value and thereby controls spatial variations of the deformation gradient without additional degrees of freedom. 
A rotation-based regularization is introduced in \cite{DahSjoDalWal:2026:ara}. 
It targets local changes of rotation, avoids penalizing stretch deformation modes and remains compatible with first-order elements without adding auxiliary fields.
In \cite{WriKorXu:2026:aic} a thin layer is introduced between the solid and the third medium. 
The idea is to make the third medium even softer with a stiffening in the thin layer. 
By this the influence of the third medium on the deformation of the solids is reduced.
Another active direction is the development of low-order-compatible formulations. 
Many gradient-based regularizations require second derivatives of the displacement field and therefore higher-order finite elements. 
This increases implementation effort and computational cost. 
First-order finite element formulations are introduced in \cite{WriKorJun:2025:atm}, where auxiliary scalar fields approximate gradients of rotation-like or skew-symmetric deformation measures. 
This enables triangular, quadrilateral, tetrahedral and hexahedral low-order elements, but introduces additional unknown fields. 
A related thermo-mechanical formulation based on low-order ansatz spaces is proposed in \cite{Wri:2026:atm}, where the third medium carries both mechanical contact forces and heat flux without explicit interface conditions. 
A stabilization-free virtual element formulation for third medium contact is presented in \cite{XuWri:2026:sfv}, allowing polygonal discretizations while avoiding the classical virtual element stabilization difficulty in the presence of higher-order regularization terms. 
A three-dimensional third medium contact framework for hyperelastic contact and pneumatically actuated systems is developed in \cite{XuXueWri:2026:tdt}. 
This formulation includes a three-dimensional regularization strategy and a pneumatic loading contribution. 
An efficient solution strategy for additional regularization fields is proposed in \cite{ZabJanJun:2026:afa}. 
There, the fields are solved in a staggered manner at quadrature points by means of a neighbored-element method, which reduces the number of global unknowns while retaining low-order displacement approximations.
The current state of third medium contact therefore provides several robust stabilization strategies. 
However, each strategy comes with specific numerical consequences. 
Some approaches require higher-order derivatives or higher-order finite elements. 
Others introduce additional global fields, rely on element-wise averaging, or use specialized discretizations. 
Low-order finite elements remain attractive since they are simple, robust and computationally cheap. 
At the same time, they require stabilization mechanisms that do not rely on the direct evaluation of second derivatives of the displacement field.

In the present work, an auxiliary-field stabilization for third medium contact formulations is introduced. 
A deformation-gradient-like field is added in the third medium and treated independently from the displacement field. 
A penalty contribution couples this auxiliary field to the physical deformation gradient, while the stabilizing term acts on the gradient of the auxiliary field rather than on the gradient of the deformation gradient itself. 
As a result, the formulation avoids the direct evaluation of second displacement derivatives and remains compatible with first- and second-order finite elements. 
Continuous and element-wise discontinuous auxiliary-field interpolations are considered in order to assess the influence of inter-element continuity on the regularization effect. 
Several benchmark problems involving large deformation, severe third-medium compression and progressive self-contact are used to evaluate the robustness and accuracy of the proposed formulation.

The remainder of this paper is organized as follows. 
First, existing third-medium contact stabilization strategies are reviewed and the present contribution is placed in context. 
The finite-deformation setting, the third-medium energy and the auxiliary-field stabilization are then introduced, followed by the finite element discretization for continuous and discontinuous auxiliary fields. 
A one-dimensional low-order example isolates the central mechanism of the method. 
Even though the deformation gradient is constant within each first-order element, a continuous auxiliary field transfers element-wise jumps into non-vanishing auxiliary-field gradients. 
This explains the effectiveness of gradient-type stabilization for low-order displacement approximations. 
The formulation is subsequently tested in two-block compression, self-contact within a box, a C-shaped self-contact problem and a three-dimensional extension. 
The paper closes with the main findings, limitations and possible extensions.

\section{Overview of existing third-medium stabilization strategies}
Several third-medium contact regularization strategies have been proposed to control excessive mesh distortion during severe compression and sliding. 
Tab.~\ref{tab:tmc_regularization_overview} summarizes the main developments relevant to the present work. 
Existing approaches either penalize higher-order displacement information directly, reduce the regularization to selected rotation-, skew- or volume-related measures, exploit special discretization technologies, or avoid additional unknowns by element-wise averaging. 
The formulation proposed here follows a different route. 
A deformation-gradient-like auxiliary field is introduced as an independently interpolated mixed variable. 
Its gradient provides the stabilizing contribution, while the penalty coupling to the physical deformation gradient controls the consistency of the auxiliary field. 
For discontinuous auxiliary interpolations, the additional degrees of freedom remain element-local and can be eliminated by static condensation.
\begin{table}[h]
\centering
\scriptsize
\setlength{\tabcolsep}{4pt}
\caption{
Overview of third-medium contact regularizations, stabilization strategies and applications.
}
\begin{tabular}{p{0.23\textwidth}|p{0.10\textwidth}|p{0.63\textwidth}}
\hline\hline
Method & Ref. & Core idea \\
\hline

Original third-medium formulation
& \cite{WriSchSch:2013:afe}
& Introduction of a fictitious intermediate medium to transfer contact forces without an explicit contact constraint or contact search. The original formulation uses a dedicated third-medium material model and establishes the basic continuum-contact idea. \\
\hline

Fictitious contact material with high-order finite elements
& \cite{BogZanKolRan:2015:ncw}
& Normal contact is represented by a fictitious contact material in combination with high-order finite elements. The approach improves the representation of smooth contact fields but does not provide a general low-order regularization mechanism. \\
\hline

Isogeometric third-medium contact
& \cite{KruNguWriDeL:2018:ifc,HuaNguZho:2018:aim}
& Extension of the third-medium concept to isogeometric and isogeometric-meshfree discretizations. Higher continuity of the approximation space helps to represent smooth deformation fields and contact transitions. \\
\hline

HuHu regularization 
& \cite{BluSigPou:2021:icm} 
& Penalization of the Hessian of the displacement field in the third medium. This controls severe element distortion during finite-strain contact but requires access to second derivatives of the displacement field. \\
\hline

HuHu-based topology optimization extensions
& \cite{FreSigPou:2024:too,BluSigPou:2023:ido}
& Use of HuHu-type void or third-medium regularization for topology optimization of self-contacting structures and nonlinear springs with internal contact. The regularization stabilizes the fictitious medium in highly compliant void regions. \\
\hline

Frictional third-medium contact
& \cite{FreRokPouSigGee:2024:aft}
& Addition of friction by an anisotropic shear contribution and a crystal-plasticity-inspired slip mechanism. The main focus is frictional force transfer rather than a new mesh-regularization measure. \\
\hline

Thermo-mechanical third-medium contact
& \cite{DalAleFrePouSig:2025:too}, \cite{Wri:2026:atm}
& Extension of third-medium contact to contact-aided thermo-mechanical regulators. Self-contact is used to tune heat transfer, while the mechanical contact treatment builds on stabilized third-medium concepts. \\
\hline

HuHu-LuLu regularization 
& \cite{FreDalSigPou:2025:itm} 
& Modification of the Hessian-based HuHu regularization by subtracting a Laplacian contribution. This reduces the penalization of bending and quadratic compression modes while retaining control of undesirable skew-type deformation modes. \\
\hline

Curvature-penalized pneumatic third medium
& \cite{FalHorDosRok:2024:tmf}
& Split of the third-medium energy into contact, regularization and pneumatic pressure contributions. A curvature penalization is used to improve the behavior of compliant third-medium regions in pneumatically actuated systems. \\
\hline

Rotation-gradient regularization 
& \cite{WriKorJun:2025:atm,DahSjoDalWal:2026:ara} 
& Penalization of gradients of rotation-related quantities, such as rotation tensors or rotation angles. The aim is to regularize excessive curvature while avoiding a direct penalization of pure stretch modes. \\
\hline

Skew-gradient regularization 
& \cite{WriKorJun:2025:fof} 
& Introduction of auxiliary fields to approximate gradients of skew-symmetric deformation measures. This enables first-order finite elements for selected regularization terms but introduces additional unknown fields. \\
\hline

Jacobian-gradient regularization 
& \cite{WriKorJun:2025:atm,DahSjoDalWal:2026:ara} 
& Penalization of spatial variations of the deformation-gradient Jacobian. The regularization targets local volumetric changes in the third medium. \\
\hline

Deformation-gradient averaging 
& \cite{FalAmaHor:2026:dga} 
& Penalization of deviations between local deformation gradients and an element-wise representative deformation gradient. Spatial variations of $\bF$ are suppressed without introducing additional degrees of freedom. \\
\hline

Stabilization-free virtual element formulation
& \cite{XuWri:2026:sfv}
& Transfer of third-medium contact to a virtual element setting. Polygonal discretizations are enabled, while classical higher-order third-medium stabilization issues are avoided by the structure of the virtual element formulation. \\
\hline

Three-dimensional third-medium framework
& \cite{XuXueWri:2026:tdt}
& Extension of third-medium contact to three-dimensional hyperelastic contact and pneumatically actuated systems. The formulation includes a three-dimensional regularization strategy and a pneumatic loading contribution. \\
\hline

Neighbored-element method 
& \cite{ZabJanJun:2026:afa} 
& Efficient treatment of additional regularization fields in a neighbored-element setting. This is primarily a solution and implementation strategy for auxiliary-field formulations rather than a new regularization measure. \\
\hline

Third medium with thin layer
& \cite{WriKorXu:2026:aic} 
& Softening of the third medium to further reduce its influence on the deformation of the contacting solids.\\
\hline

Present approach 
& present work
& Introduction of a deformation-gradient-like auxiliary field $\BTheta$ that is coupled directly to the physical deformation gradient $\bF$ via a penalty term. The gradient $\nabla\BTheta$ provides the regularization measure. Continuous variants introduce inter-element stabilization, while discontinuous variants allow element-level static condensation. \\
\hline\hline
\end{tabular}
\label{tab:tmc_regularization_overview}
\end{table}

%=========================================================================
\section{Continuum mechanical foundations}
\vspace{-4mm}
%=========================================================================

Within this section, the continuum formulation of the proposed third medium contact approach is introduced. 
The physical solid domain is denoted by $\B_{\mathrm{s}}$, while the space between potentially contacting surfaces is represented by the fictitious third medium domain $\B_{\mathrm{tm}}$. 
Both domains are described within finite deformation kinematics. 
For the solid, a standard hyperelastic material model is employed. 
In the third medium, the same hyperelastic base energy is used in a strongly scaled form and supplemented by an additional stabilization contribution.

\subsection{Boundary value problem, kinematics and weak form}
All balance equations are formulated with respect to the reference configuration $\B \subset \mathbb{R}^3$. 
Under static loading, the balance of linear momentum reads
\begin{equation}
	\Div{\bP}
	+
	\rho_0 \bb
	=
	\bzero
	\qquad \mathrm{in}\quad \B ,
	\label{eq:BaMo}
\end{equation}
where $\bP$ is the first Piola-Kirchhoff stress tensor, $\rho_0$ is the mass density in the reference configuration and $\bb$ denotes the body force per unit mass. 
For the boundary, the standard decomposition
\begin{equation}
	\partial\B
	=
	\partial\B_{u}
	\cup
	\partial\B_{t},
	\qquad
	\partial\B_{u}
	\cap
	\partial\B_{t}
	=
	\emptyset
\end{equation}
is used, together with
\begin{equation}
	\bu=\bar{\bu}
	\quad\mathrm{on}\quad
	\partial\B_{u}
	\qquad\textrm{and}\qquad
	\bP\cdot\bN=\bar{\bt}
	\quad\mathrm{on}\quad
	\partial\B_{t}.
	\label{eq:BC}
\end{equation}
Here, $\bar{\bu}$ is the prescribed displacement, $\bar{\bt}$ is the prescribed traction and $\bN$ denotes the outward unit normal in the reference configuration.
A material point $\bX\in\B$ is mapped to the current configuration by
\begin{equation}
	\bx
	=
	\Bvarphi(\bX,t)
	=
	\bX+\bu(\bX,t).
\end{equation}
Consequently, the deformation gradient, its Jacobian and the right Cauchy-Green tensor are given by
\begin{equation}
	\bF
	=
	\frac{\partial \bx}{\partial \bX}
	=
	\bI+\nabla\bu ,
	\qquad
	J
	=
	\det\bF 
	\qquad\textrm{and}\qquad
	\bC
	=
	\bF^{T}\cdot\bF .
	\label{eq:F}
\end{equation}
$J$ measures the local volume change. 
The volume-preserving part of the deformation is described by the isochoric deformation gradient
\begin{equation}
	\widehat{\bF}
	=
	J^{-1/3}\bF .
\end{equation}
%
%In the solid domain, a compressible Neo-Hookean strain energy density is used as
In the solid domain, the strain-energy density is denoted by $\psi$. 
Here, a compressible Neo-Hookean model is used as
\begin{equation}
	\psi(\bC)
	=
	\frac{K}{2}
	\left(\ln J\right)^2
	+
	\frac{\mu}{2}
	\left(
	J^{-2/3}\tr\bC
	-
	3
	\right),
	\label{eq:NeoHook}
\end{equation}
where $K$ and $\mu$ are the bulk and shear moduli, respectively. 
Stress measures follow from
\begin{equation}
	\bS
	=
	2
	\frac{\partial\psi(\bC)}{\partial\bC},
	\qquad
	\bP
	=
	\bF\cdot\bS .
	\label{eq:PK}
\end{equation}
Multiplication of Eq.~\ref{eq:BaMo} with an admissible virtual displacement $\delta\bu$ and integration over the reference domain yields the weak form
\begin{equation}
\begin{aligned}
	\delta\Pi
	=
	\int_{\B}
	\frac{\partial\psi(\bC)}{\partial\bF}
	:
	\delta\bF
	\,\mathrm{d}V
	-
	\int_{\B}
	\rho_0
	\bb
	\cdot
	\delta\bu
	\,\mathrm{d}V
	-
	\int_{\partial\B_t}
	\bar{\bt}
	\cdot
	\delta\bu
	\,\mathrm{d}A
	=
	0 .
\end{aligned}
\label{eq:dynamic_weak_form}
\end{equation}
Here, $\delta\bF=\nabla\delta\bu$ denotes the variation of the deformation gradient. 
Body forces $\bb$ are neglected within the following.

\subsection{Third medium energy}

In the third medium, the corresponding strain-energy density $\psi$ is used as a fictitious base material that is scaled by the small parameter $\gamma$, while $W^{\mathrm{tm}}$ denotes the total third-medium energy obtained by integration over $\B_{\mathrm{tm}}$.
This scaling keeps the influence of the third medium negligible before contact, whereas the underlying energy still increases strongly under severe compression and thereby provides the desired contact-barrier effect.
For three-dimensional simulations, the third-medium energy is defined as
\begin{equation}
	W^{\mathrm{tm}}
	=
	\gamma
	\int_{\B_{\mathrm{tm}}}
	\left[
	\frac{K}{2}
	\left(\ln J\right)^2
	+
	\frac{\mu}{2}
	\left(
	J^{-2/3}
	\tr\bC
	-
	3
	\right)
	\right]
	\,\mathrm{d}V .
	\label{eq:NeoHookTM}
\end{equation}
Volumetric changes are controlled by the first term, while the second term describes the isochoric deformation of the third medium. 
For the three-dimensional computations, the full third-medium energy in Eq.~\ref{eq:NeoHookTM} is retained. 
In the two-dimensional plane strain examples, however, a reduced contact-barrier energy is used. 
This choice follows the common interpretation of the third medium as a fictitious contact medium rather than as a physical material phase: its stiffness should remain negligible before contact, while the energy has to increase strongly when the intermediate region is compressed towards a nearly collapsed state, cf.~\cite{WriSchSch:2013:afe,WriKorJun:2025:atm}. 
Within the present plane strain setting, the deformation is embedded in three dimensions by setting \(F_{33}=1\). 
The isochoric contribution then already provides the required barrier effect for strongly collapsing in-plane deformations, since it becomes singular as \(J \to 0\). 
Consequently, the volumetric contribution is omitted in the reduced two-dimensional model, leading to
\begin{equation}
W^{\mathrm{tm}}_{\mathrm{red}}
=
\gamma
\int_{\mathcal{B}_{\mathrm{tm}}}
\frac{\mu}{2}
\left(
J^{-2/3}\operatorname{tr}\mathbf{C}-3
\right)
\,\mathrm{d}V .
\label{eq:Wtm_red}
\end{equation}
Similar reduced two-dimensional third-medium energies based on distortional or isochoric contributions have been used in \cite{ZabJanJun:2026:afa}. 
It should be emphasized that Eq.~\ref{eq:Wtm_red} is not intended as a complete three-dimensional material law. 
It is used only as a reduced two-dimensional contact-barrier energy in the plane strain examples, whereas all three-dimensional simulations employ the full energy in Eq.~\ref{eq:NeoHookTM}.

\subsection{Deformation-gradient-based stabilization}
A useful way to motivate the proposed stabilization is to consider a gradient-type control of the deformation gradient in the third medium. 
Such a term would penalize spatial variations of $\bF$ and can formally be written as

\begin{equation}
	W^{\mathrm{tm}}_{\nabla\bF}
	=
	\int_{\B_{\mathrm{tm}}}
	\frac{\alpha_r}{2}
	\left|
	\nabla \bF
	\right|^2
	\,\mathrm{d}V ,
	\label{eq:TMC_direct_gradF}
\end{equation}
where $\alpha_r$ is a regularization parameter. 
Here, Eq.~\ref{eq:TMC_direct_gradF} is used only as a motivating reference form and is not meant to represent the specific structure of all existing third-medium regularizations. 
Instead, it highlights the quantity that the present formulation aims to control without evaluating second derivatives of the displacement field directly.
For low-order finite elements, however, such a direct term is not suitable. 
With a first-order displacement interpolation, the deformation gradient is element-wise constant and its gradient cannot provide an effective intra-element stabilization. 
To avoid the direct evaluation of $\nabla \mathbf{F}$, an additional deformation-gradient-like field $\BTheta$ is introduced in the third medium. 
This field is interpolated independently from the displacement field and is weakly coupled to the physical deformation gradient by the penalty contribution as 
\begin{equation}
	W^{\mathrm{tm}}_{p}
	=
	\int_{\B_{\mathrm{tm}}}
	\frac{p_\Theta}{2}
	\left|
	\BTheta-\bF
	\right|^2
	\,\mathrm{d}V ,
	\label{eq:TMC_penalty_theta_F}
\end{equation}
where $p_\Theta$ is the penalty parameter. 
For sufficiently large $p_\Theta$, the auxiliary field is driven towards the deformation gradient. 
Regularization is then applied to the spatial gradient of the auxiliary field,
\begin{equation}
	W^{\mathrm{tm}}_{r}
	=
%	\gamma
	\int_{\B_{\mathrm{tm}}}
	\frac{\alpha_r}{2}
	\left|
	\nabla\BTheta
	\right|^2
	\,\mathrm{d}V .
	\label{eq:TMC_grad_theta}
\end{equation}
In this way, a gradient regularization of $\bF$ is approximated without computing $\nabla\bF$ directly. 
Only first derivatives of the independently interpolated field $\BTheta$ are required. 
For three-dimensional simulations, the stabilized third-medium energy becomes
\begin{equation}
\begin{aligned}
	W^{\mathrm{tm}}_{\mathrm{stab}}
	&=
	W^{\mathrm{tm}}
	+
	W^{\mathrm{tm}}_{p}
	+
	W^{\mathrm{tm}}_{r}
	\\[2mm]
	&=
	\int_{\B_{\mathrm{tm}}}
	\bigg[
	\gamma
	\left[
	\frac{K}{2}
	\left(\ln J\right)^2
	+
	\frac{\mu}{2}
	\left(
	J^{-2/3}
	\tr\bC
	-
	3
	\right)
	\right]
	+
	\frac{p_\Theta}{2}
\norm{\BTheta-\bF}^2
	+
%	\gamma
	\frac{\alpha_r}{2}
    \norm{\nabla\BTheta}^2
	\bigg]
	\,\mathrm{d}V .
\end{aligned}
\label{eq:TMC_total_energy_theta}
\end{equation}
For the reduced two-dimensional setting, $W^{\mathrm{tm}}$ in Eq.~\ref{eq:TMC_total_energy_theta} is replaced by $W^{\mathrm{tm}}_{\mathrm{red}}$ from Eq.~\ref{eq:Wtm_red}. 
Hence, the proposed stabilization provides a low-order-compatible approximation of a deformation-gradient regularization while avoiding the direct computation of second derivatives of the displacement field.

\subsection{Finite element discretization and static condensation}
The proposed stabilization leads to a mixed finite element formulation in the third medium. 
The displacement field $\bu$ and the deformation-gradient-like auxiliary field $\BTheta$ are treated as independent unknowns. 
The physical deformation gradient $\bF$ remains kinematically tied to the displacement field, while $\BTheta$ serves as an additional stabilization field. 
Both fields are coupled weakly through the penalty contribution in Eq.~\ref{eq:TMC_penalty_theta_F}. 
The regularization acts on $\nabla\BTheta$ and therefore avoids the direct use of $\nabla\bF$. 
Consequently, the formulation mimics a deformation-gradient regularization without requiring second derivatives of the displacement field.
Within an isoparametric finite element setting, the displacement field and the auxiliary field are interpolated independently. 
On element level, the approximations read
\begin{equation}
	\underline{\bu}
	=
	\underline{\IN}^{e}
	\underline{\bd}_{u}^{e},
	\qquad
	\underline{\BTheta}
	=
	\underline{\IN}_{\Theta}^{e}
	\underline{\bd}_{\Theta}^{e}
	\label{eq:FE_interpolation}
\end{equation}
while associated gradient quantities are computed as
\begin{equation}
	\underline{\bF}
	=
	\underline{\bI}
	+
	\underline{\IB}_{u}^{e}
	\underline{\bd}_{u}^{e},
	\qquad
	\nabla\underline{\BTheta}
	=
	\underline{\IB}_{\Theta}^{e}
	\underline{\bd}_{\Theta}^{e}.
	\label{eq:FE_gradients}
\end{equation}
A central feature of the formulation is that the interpolation order of $\BTheta$ is not tied to the interpolation order of $\bu$. 
The auxiliary field may therefore be approximated with the same polynomial degree as the displacement field, but lower-order approximations are equally possible. 
Since $\BTheta$ serves as a numerical stabilization field rather than as an independent physical field, no physical continuity requirement is imposed a priori. 
This allows both continuous and element-wise discontinuous auxiliary interpolations to be considered within the same framework.
In the present work, the interpolation pairs 
$\mathrm{T}_{1}^{u}\mathrm{T}_{1}^{\Theta}$, 
$\mathrm{T}_{2}^{u}\mathrm{T}_{2}^{\Theta}$ and 
$\mathrm{T}_{2}^{u}\mathrm{T}_{1}^{\Theta}$ 
are compared with their discontinuous counterparts
$\mathrm{T}_{1}^{u}\mathrm{T}_{1}^{\Theta,d}$, 
$\mathrm{T}_{1}^{u}\mathrm{T}_{0}^{\Theta,d}$, and 
$\mathrm{T}_{2}^{u}\mathrm{T}_{1}^{\Theta,d}$. 
Here, the subscript denotes the polynomial degree, while the superscript identifies the approximated field. 
No additional superscript indicates a continuous auxiliary interpolation, whereas the superscript $d$ denotes an element-wise discontinuous approximation of $\BTheta$.
For discontinuous auxiliary interpolations, all $\BTheta$-degrees of freedom remain local to the element. 
Static condensation can therefore be performed on element level, such that the global Newton system contains only displacement degrees of freedom. 
For continuous auxiliary interpolations, the corresponding degrees of freedom are globally coupled and are therefore assembled in the standard way.
After linearization, the element contribution takes the block form
\begin{equation}
\begin{bmatrix}
	\underline{\IK}^{e}_{uu}
	&
	\underline{\IK}^{e}_{u\Theta}
	\\[2mm]
	\underline{\IK}^{e}_{\Theta u}
	&
	\underline{\IK}^{e}_{\Theta\Theta}
\end{bmatrix}
\begin{bmatrix}
	\Delta\underline{\bd}^{e}_{u}
	\\[2mm]
	\Delta\underline{\bd}^{e}_{\Theta}
\end{bmatrix}
=
-
\begin{bmatrix}
	\underline{\IR}^{e}_{u}
	\\[2mm]
	\underline{\IR}^{e}_{\Theta}
\end{bmatrix}.
\label{eq:block_system}
\end{equation}

For continuous auxiliary-field interpolations, both fields contribute to the global system, and Eq.~\ref{eq:block_system} represents the element-level form of the coupled mixed problem. 
For discontinuous auxiliary-field interpolations, all auxiliary degrees of freedom are element-local. 
In this case, the mixed character is retained at element level, but the auxiliary unknowns can be eliminated before global assembly.
For element-wise discontinuous auxiliary fields, elimination of the auxiliary increment yields the condensed displacement equation
\begin{equation}
	\underline{\IK}^{e,\mathrm{c}}_{uu}
	\,
	\Delta\underline{\bd}^{e}_{u}
	=
	-
	\underline{\IR}^{e,\mathrm{c}}_{u},
	\label{eq:condensed_system}
\end{equation}
with
\begin{equation}
	\underline{\IK}^{e,\mathrm{c}}_{uu}
	=
	\underline{\IK}^{e}_{uu}
	-
	\underline{\IK}^{e}_{u\Theta}
	\left(
	\underline{\IK}^{e}_{\Theta\Theta}
	\right)^{-1}
	\underline{\IK}^{e}_{\Theta u}
\quad\textrm{and}\quad	
	\underline{\IR}^{e,\mathrm{c}}_{u}
	=
	\underline{\IR}^{e}_{u}
	-
	\underline{\IK}^{e}_{u\Theta}
	\left(
	\underline{\IK}^{e}_{\Theta\Theta}
	\right)^{-1}
	\underline{\IR}^{e}_{\Theta}.
\label{eq:schur_complement}
\end{equation}
Once the global displacement increment is known, the auxiliary increment is recovered locally by
\begin{equation}
	\Delta\underline{\bd}^{e}_{\Theta}
	=
	-
	\left(
	\underline{\IK}^{e}_{\Theta\Theta}
	\right)^{-1}
	\left(
	\underline{\IR}^{e}_{\Theta}
	+
	\underline{\IK}^{e}_{\Theta u}
	\Delta\underline{\bd}^{e}_{u}
	\right)
	\label{eq:theta_increment_recovery}
\end{equation}
and the auxiliary degrees of freedom are updated according to
\begin{equation}
	\underline{\bd}^{e,it+1}_{\Theta}
	=
	\underline{\bd}^{e,it}_{\Theta}
	+
	\Delta\underline{\bd}^{e,it}_{\Theta}.
	\label{eq:DOF_update}
\end{equation}
Thus, the discontinuous variants retain a local auxiliary-field coupling without increasing the size of the global system, while the continuous variants provide the stronger inter-element regularization mechanism investigated below.

%=========================================================================
\section{Numerical examples}
\vspace{-4mm}
%=========================================================================
%
Within this section the performance of the proposed auxiliary-field stabilization is assessed. 
The numerical study is organized around benchmark problems that test different aspects of the formulation, including the influence of the third-medium scaling, the penalty coupling, the gradient regularization and the interpolation order. 
First, two elastic blocks are pressed into contact to analyze the role of the third-medium parameters and the residual gap. 
Further verification is provided by the self-contact-within-a-box benchmark and the classical C-shaped boundary value problem, for which reference results are available in the literature. 
Most examples are formulated in two dimensions in order to isolate the relevant stabilization effects and allow a clear comparison of the different interpolation pairs. 
A three-dimensional version of the self-contact-within-a-box problem is included to demonstrate that the formulation is not restricted to plane settings.

\subsection{$\nabla\BTheta$-regularization for continuous low-order interpolation \label{sec:1DRegularization}}
\begin{wrapfigure}{r}{8cm}
%\begin{Figure}[H]
\begin{picture}(0,6)
\unitlength1cm
\put( 0.5,3.8){\includegraphics[width=8cm]{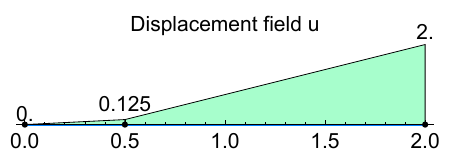}}
\put( 0.0,4.0){a)}
\put( 0.5,1.0){\includegraphics[width=8cm]{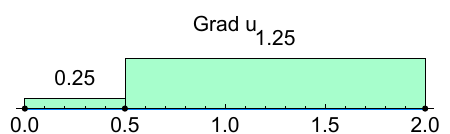}}
\put( 0.0,1.2){b)}
\put( 0.5,-2.5){\includegraphics[width=8cm]{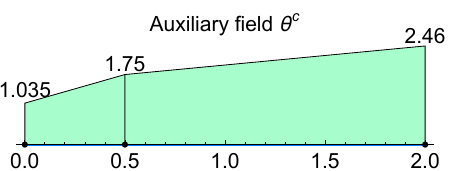}}
\put( 0.0,-2.){c)}
\put( 0.5,-6.2){\includegraphics[width=8cm]{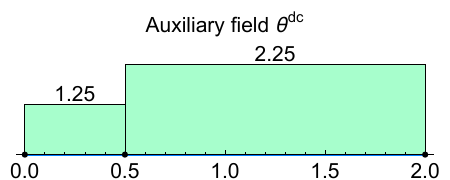}}
\put( 0.0,-5.8){d)}
\end{picture}
\vspace*{63mm}
\caption{
One-dimensional illustration of the auxiliary-field regularization for two linear finite elements. 
a) Prescribed displacement field,
b) gradient of u,
c) auxiliary field with a continuous interpolation, and
d) auxiliary field with a discontinuous interpolation.
}
\label{fig:1DGradient}
%\end{Figure}
\end{wrapfigure}
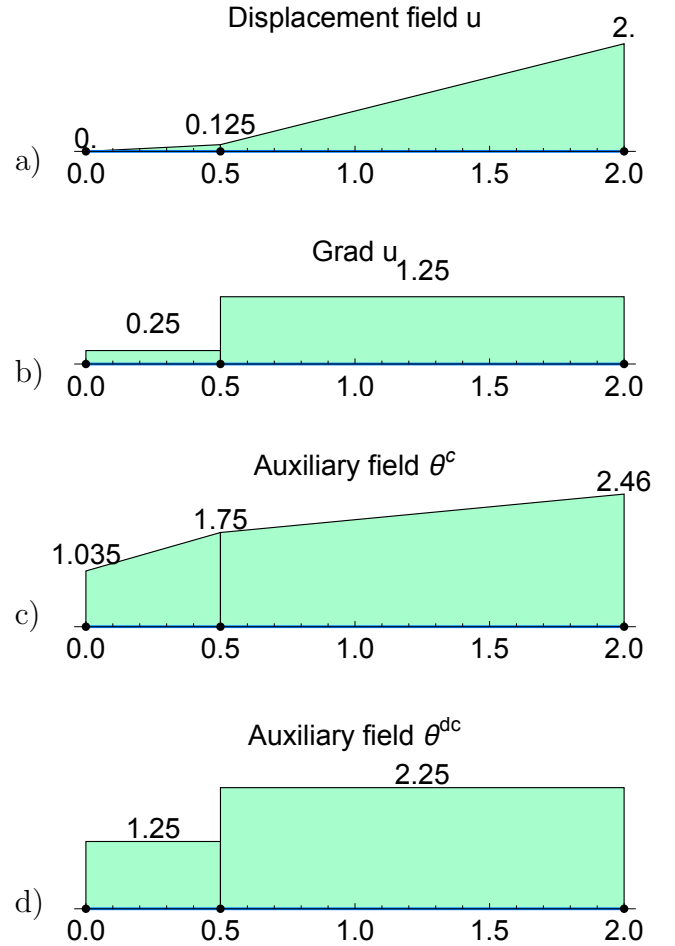

This example illustrates why a low-order discretization can still generate a non-vanishing gradient of the auxiliary field $\BTheta$, and why this mechanism requires a continuous interpolation of $\BTheta$. 
For clarity, the setting is reduced to a one-dimensional boundary value problem discretized by two linear finite bar elements. 
The total length is $L=2$, with element lengths $L_1=0.25$ and $L_2=1.75$. 
Three geometrical nodes are obtained, and the displacement field is fully prescribed by $u^1=0$, $u^2=1/8$ and $u^3=7/4$. 
Hence, the deformation field is fixed and only the auxiliary degrees of freedom associated with $\BTheta$ remain unknown. 
Depending on the interpolation chosen for $\BTheta$, different field distributions arise inside the elements and across the shared element boundary.
Fig.~\ref{fig:1DGradient} shows the prescribed displacement field and the resulting auxiliary-field distributions. 
Due to the linear displacement interpolation, the deformation gradient is constant within each element, but differs between the short left element and the longer right element. 
With a continuous interpolation of $\BTheta$, the auxiliary field shares one degree of freedom at the common node. 
This shared value has to represent the deformation-gradient information of both adjacent elements and therefore acts as an inter-element transition value. 
As a result, $\BTheta$ becomes piecewise linear and develops a non-zero gradient across the two elements, as shown in Fig.~\ref{fig:1DGradient}b).
For an element-wise discontinuous interpolation, each element owns its local auxiliary degrees of freedom. 
The auxiliary field can then represent the constant deformation gradient of each element independently, without enforcing compatibility at the shared geometrical node. 
Consequently, the field becomes discontinuous at the element interface, as shown in Fig.~\ref{fig:1DGradient}c), and no continuity-driven gradient-based coupling between neighboring elements is generated. 
Thus, the continuous interpolation introduces a non-local stabilization mechanism: changes of the deformation gradient between adjacent elements are transferred into gradients of $\BTheta$ and can be penalized by the regularization term in Eq.~\ref{eq:TMC_grad_theta}. 
This observation is relevant not only for the present formulation, but also for related third-medium contact stabilizations that rely on auxiliary fields or gradient-type regularization measures.

\subsection{Self-contact between two blocks}
In this example the influence of the third-medium scaling parameter $\gamma$ on the contact response is investigated. 
Two elastic blocks of size $100\times50$ are separated by a third medium layer of the same in-plane size, as shown in Fig.~\ref{fig:BVPcontactBlocksRef}. 
Dirichlet boundary conditions are prescribed on the top and bottom boundaries. 
The lower boundary is fixed by $\overline{\bu}=\bzero$, while the upper boundary is displaced according to two loading cases. 
A T$_1^u$T$_1^\Theta$ discretization is applied. 

\begin{Figure}[h]
\begin{picture}(0,7.5)
\unitlength1cm
\put( 3.0,1.0){\includegraphics[width=0.3\textwidth]{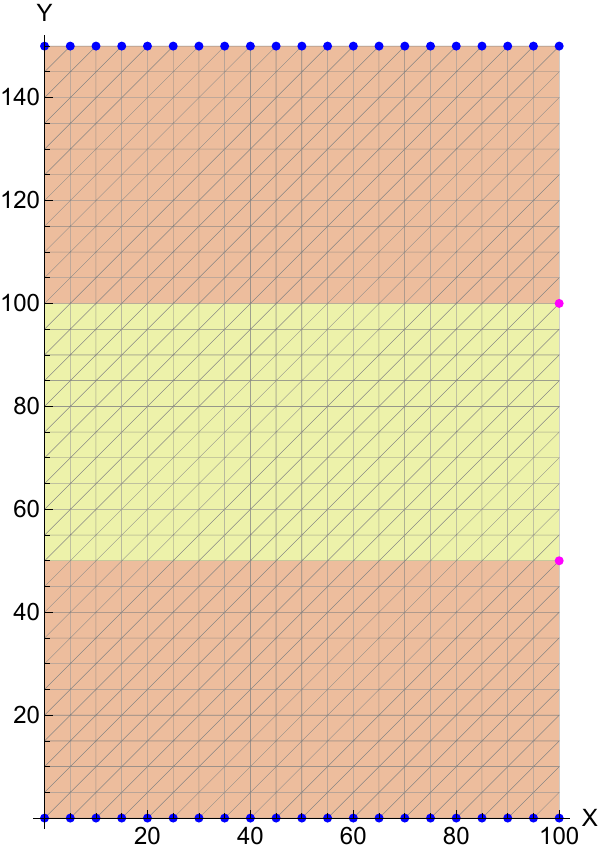}}
\put(10.0,6.0){\textbf{Parameters}}
\put(10.3,5.3){Shear modulus: $\mu=5/14$}
\put(10.3,4.8){Bulk modulus: $K=5/3$}
\put(10.0,3.8){\textbf{TMC parameters}}
\put(10.3,3.1){$\gamma=10^{-3},10^{-4},10^{-5},10^{-8}$}
\put(10.3,2.6){$\alpha_r=1$}
\put(10.3,2.1){$p_\Theta=1$}
\put(7.6,5.3){$u_2^o$}
\put(7.6,3.3){$u_2^u$}
\end{picture}
\vspace*{-8mm}
\caption{
Reference configuration of the two-block contact benchmark. 
The upper and lower solid blocks are separated by a third medium layer, which transfers the contact forces during compression. 
Dirichlet boundary conditions are prescribed on the top and bottom boundaries. 
The lower boundary is fixed with $\overline{\bu}=\bzero$, while the upper boundary is displaced in case 1 until $\overline{\bu}=[0,-60]^T$ and in case 2 until $\overline{\bu}=[10,-60]^T$. 
$u_2^o$ and $u_2^u$ indicate the $u_2$-displacements within the marked nodes. 
}
\label{fig:BVPcontactBlocksRef}
\end{Figure}

In case~1, a purely vertical compression is applied with $\overline{\bu}=[0,-60]^T$. 
In case~2, a combined horizontal and vertical compression is applied with $\overline{\bu}=[10,-60]^T$.
The solid blocks are modeled by the compressible Neo-Hookean material in Eq.~\ref{eq:NeoHook} with $K_{\mathrm{s}}=5/3$ and $\mu_{\mathrm{s}}=5/14$. 
Inside the third medium, the same base energy is used with $\gamma=\{10^{-3},10^{-4},10^{-5},10^{-8}\}$ and the stabilization parameters $p_{\Theta}=1$ and $\alpha_r=1$. 
This setup is used to assess how the stiffness scaling of the third medium affects the residual gap between both blocks.
For the evaluation, two reference nodes are tracked at the right edge of the third medium. 
The quantity $u_2^o$ denotes the vertical displacement of the upper reference node $P(100,100)$, while $u_2^u$ denotes the vertical displacement of the lower reference node $P(100,50)$. 
The difference $u_2^o-u_2^u$ is used as a scalar measure for the remaining gap between the two blocks. 

%\newpage

\begin{Figure}[h]
\begin{picture}(0,6.3)
\unitlength1cm
\put( 1.0,1.0){\includegraphics[width=0.4\textwidth]{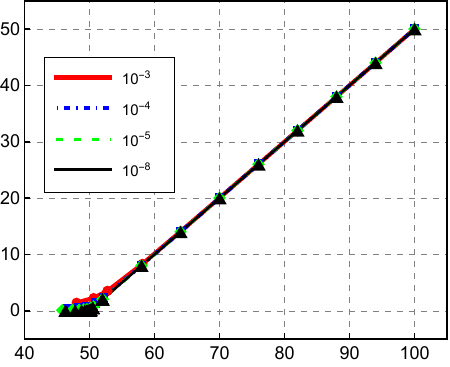}}
\put( 1.0,0.6){a)}
\put( 8.5,1.0){\includegraphics[width=0.4\textwidth]{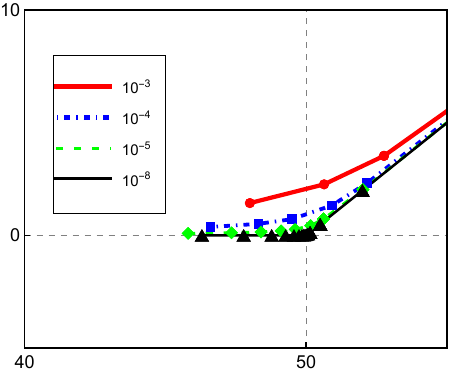}}
\put( 8.5,0.6){b)}
\put( 3.9,0.6){$u_2^o-u_2^u$}
\put( 0.5,3.5){\rotatebox{90}{$u_2$}}
\put(11.1,0.6){$u_2^o-u_2^u$}
\put( 8.3,3.5){\rotatebox{90}{$u_2$}}
\end{picture}
\vspace*{-6mm}
\caption{
Displacement-gap curves for the two-block benchmark under vertical compression. 
The gap measure is defined by $u_2^o-u_2^u$, where $u_2^o$ and $u_2^u$ denote the vertical displacements of the marked reference nodes. 
Results are shown for different values of the third-medium scaling parameter $\gamma$ given in a) as the full loading path and b) as a magnified perspective of the contact regime.
}
\label{fig:contactPlotVertical}
\end{Figure}

Fig.~\ref{fig:contactPlotVertical} reports the displacement-gap response for case~1 and different values of $\gamma$. 
At the beginning of the loading path, the curves are almost identical, which confirms the weak influence of the third medium before contact. 
After the gap is closed, the curves separate and the effect of the scaling parameter becomes visible.

\begin{Figure}[h]
\begin{picture}(0,6)
\unitlength1cm
\put( 0.0,1.0){\includegraphics[width=0.24\textwidth]{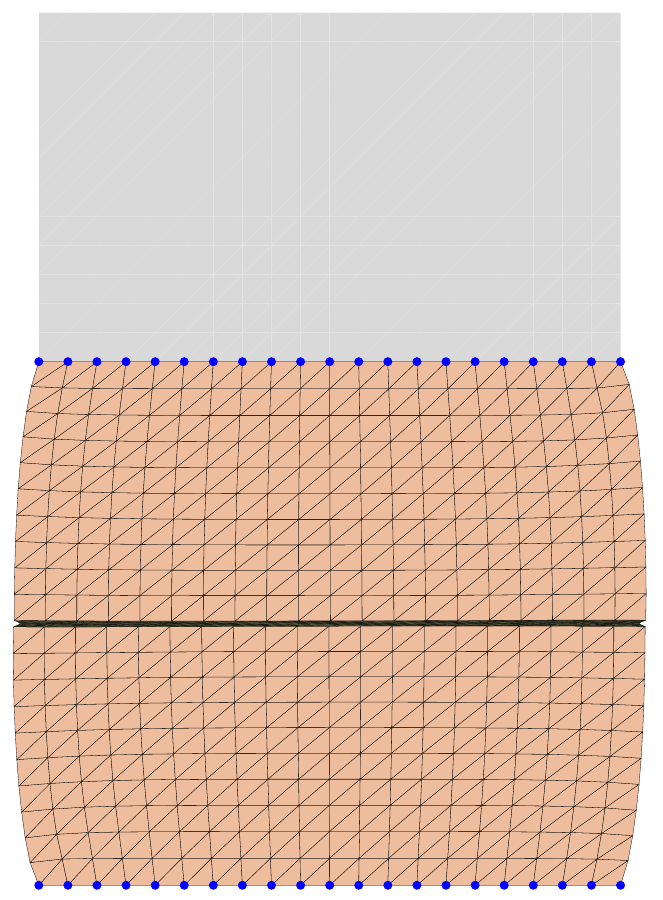}}
\put( 0.0,0.5){a) $\gamma=10^{-3}$}
\put( 4.,1.0){\includegraphics[width=0.24\textwidth]{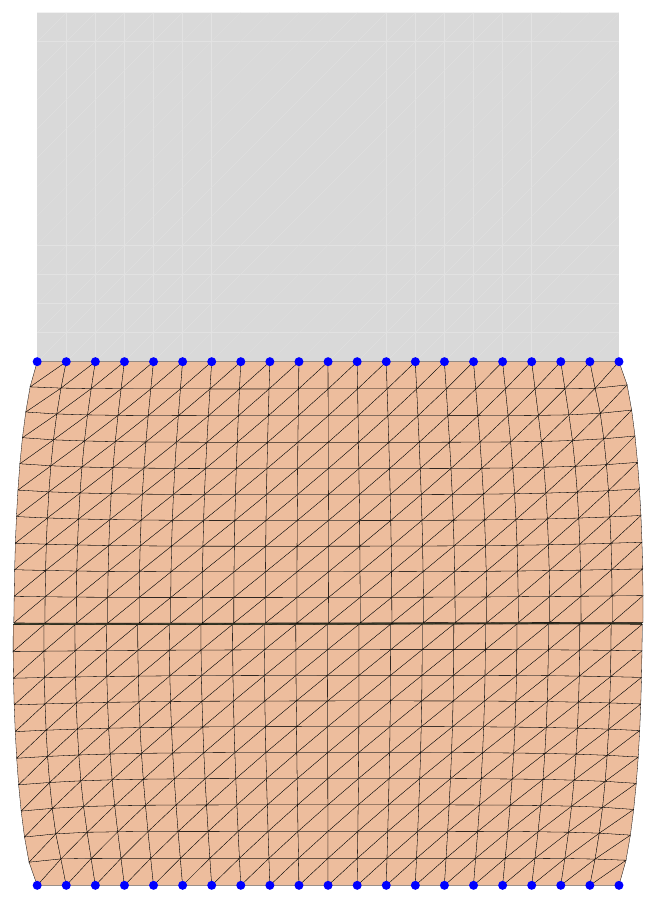}}
\put( 4.0,0.5){b) $\gamma=10^{-4}$}
\put( 8,1.0){\includegraphics[width=0.24\textwidth]{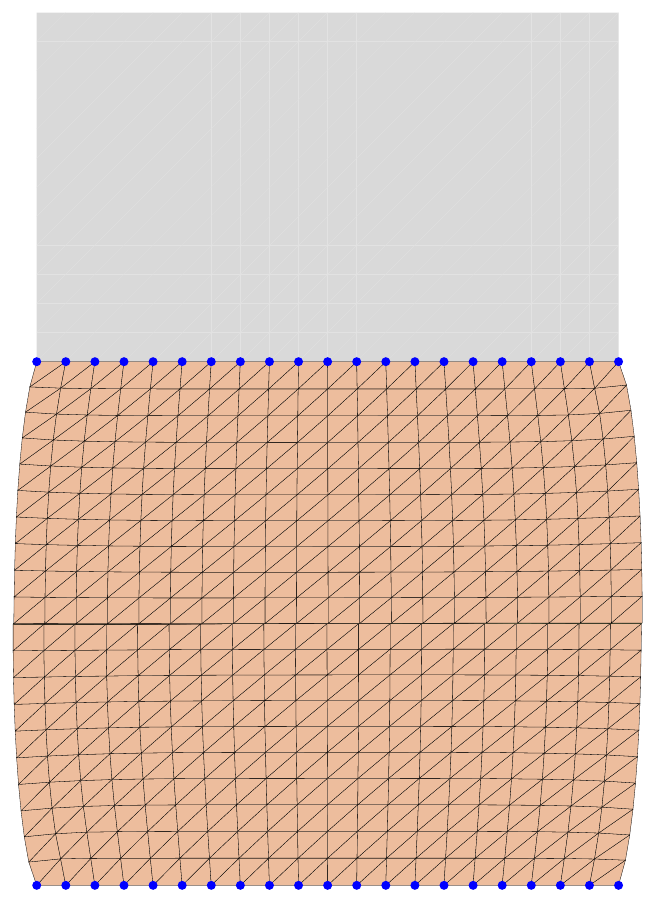}}
\put( 8.0,0.5){c) $\gamma=10^{-5}$}
\put(12,1.0){\includegraphics[width=0.24\textwidth]{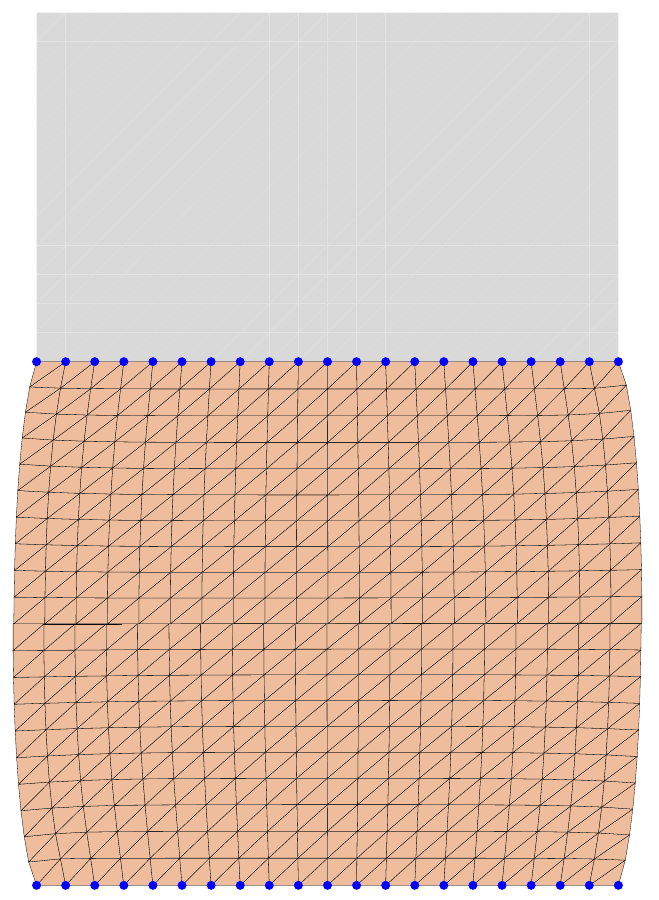}}
\put(12.0,0.5){d) $\gamma=10^{-8}$}
\end{picture}
\vspace*{-6mm}
\caption{
Final deformed configurations of the two-block benchmark under vertical compression for different values of the third-medium scaling parameter $\gamma$. 
The gray overlay indicates the reference configuration.
}
\label{fig:contactBlocksVertDeformed}
\end{Figure}

The trend observed in Fig.~\ref{fig:contactPlotVertical} can be interpreted directly from the role of the scaling parameter $\gamma$. 
This parameter controls the stiffness level of the fictitious third medium relative to the surrounding solid. 
Larger values of $\gamma$ make the intermediate layer more resistant against compression and therefore lead to a more pronounced residual gap after contact has been established, see also Fig.~\ref{fig:contactBlocksVertDeformed}.
Reducing $\gamma$ weakens the third medium before contact and allows the two solid bodies to approach each other more closely, which results in a smaller residual gap. 
This improved contact closure, however, comes at the price of a more demanding nonlinear solution process, since the third medium has to sustain severe compression while carrying only a very small stiffness contribution away from contact. 
The increasing number of Newton iterations in Tab.~\ref{tab:convstudy} reflects this trade-off between contact accuracy and nonlinear robustness. 
Within the investigated range, all converged curves remain smooth in the contact regime, indicating a stable transition from free motion to contact force transmission.

\begin{Table}[h!]
\centering
\caption{
Residual gap and total Newton iteration count for the two-block benchmark under vertical compression. 
Results are listed for different combinations of the third-medium scaling parameter $\gamma$ and the penalty parameter $p_\Theta$. 
A dash denotes a non-converged parameter set.
}
\begin{tabular}{c|c|c|c|||c|c|c|c}
  \hline\hline
  \rule{0pt}{13pt}\hspace*{-1mm}
$p_\Theta$&$\gamma$ & gap  & iter. &$p_\Theta$&$\gamma$ & gap  & iter. \\[1mm]
  \hline\hline\rule{0pt}{13pt}\hspace*{-1mm}
1&$10^{-3}$ & 1.4377 & 54 &$10^{-1}$&$10^{-3}$ & - &  - \\[1mm]
%  \hline\rule{0pt}{13pt}\hspace*{-1mm}
&$10^{-4}$ & 0.3748 &  66 &&$10^{-4}$ & 0.2846 & 68\\[1mm]
%  \hline\rule{0pt}{13pt}\hspace*{-1mm}
&$10^{-5}$&0.0908&   82 &&$10^{-5}$ & 0.0847 & 82 \\[1mm]  
%  \hline\rule{0pt}{13pt}\hspace*{-1mm}
&$10^{-8}$& 0.0015 &  115 &&$10^{-8}$ & 0.0014 & 118\\[1mm]
  \hline\hline
\end{tabular}
\label{tab:convstudy}
\end{Table}

This visual impression is quantified in Tab.~\ref{tab:convstudy}, where the final residual gap and the total number of Newton iterations are listed for different penalty parameter combinations of $\gamma$ and $p_\Theta$.
For $p_\Theta=1$, reducing $\gamma$ from $10^{-3}$ to $10^{-8}$ decreases the residual gap from $1.4377$ to $0.0015$. 
The same trend is observed for $p_\Theta=10^{-1}$, where the gap values remain close to those obtained with $p_\Theta=1$ for the converged cases. 
The number of Newton iterations increases as $\gamma$ decreases, which reflects the stronger enforcement of the contact-like response through a softer third medium. 
For the largest value $\gamma=10^{-3}$ and the smaller penalty $p_\Theta=10^{-1}$, no converged solution is obtained.

\begin{Figure}[h]
\begin{picture}(0,6.5)
\unitlength1cm
\put( 1.0,1.0){\includegraphics[width=0.4\textwidth]{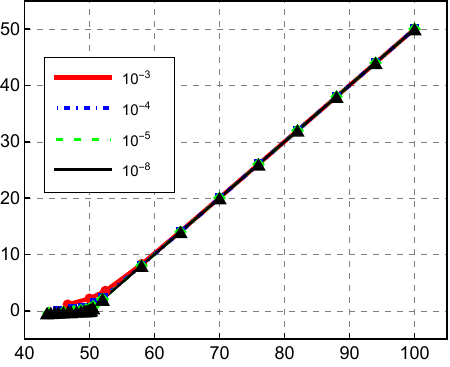}}
\put( 1.0,0.6){a)}
\put( 8.5,1.0){\includegraphics[width=0.4\textwidth]{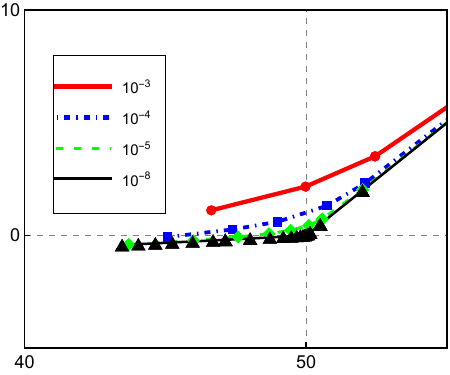}}
\put( 8.5,0.6){b)}
\put( 3.9,0.6){$u_2^o-u_2^u$}
\put( 0.5,3.5){\rotatebox{90}{$u_2^o$}}
\put(11.1,0.6){$u_2^o-u_2^u$}
\put( 8.3,3.5){\rotatebox{90}{$u_2^o$}}
\end{picture}
\vspace*{-6mm}
\caption{
Displacement-gap curves for the two-block benchmark under combined horizontal and vertical loading. 
The gap measure is defined by $u_2^o-u_2^u$, where $u_2^o$ and $u_2^u$ denote the vertical displacements of the marked reference nodes in Fig.~\ref{fig:BVPcontactBlocksRef}. 
Results are shown for different values of the third-medium scaling parameter $\gamma$. 
(a) Full loading path. 
(b) Magnified view of the contact regime.
}
\label{fig:contactPlotShear}
\end{Figure}

For case~2, the upper block is subjected to a combined horizontal and vertical displacement with $\overline{\bu}=[10,-60]^T$. 
Compared with the purely vertical loading in case~1, this load case additionally introduces tangential motion between the two blocks and therefore activates a more demanding shear-dominated deformation state in the third medium. 
Fig.~\ref{fig:contactPlotShear} shows the corresponding displacement-gap response for different values of $\gamma$. 
The overall trend remains consistent with case~1: before contact, the curves are close to each other, while the influence of $\gamma$ becomes visible once the gap is closed. 
In the contact regime, larger values of $\gamma$ again lead to a larger residual gap, whereas smaller values allow a tighter closure of the interface. 
The zoom in Fig.~\ref{fig:contactPlotShear}b highlights that the response remains smooth even under combined compression and tangential loading. 
This indicates that the proposed stabilization remains robust when the third medium is exposed not only to compression, but also to pronounced shear deformation.

\begin{Figure}[h]
\begin{picture}(0,5.8)
\unitlength1cm
\put( 0.0,1.0){\includegraphics[width=0.24\textwidth]{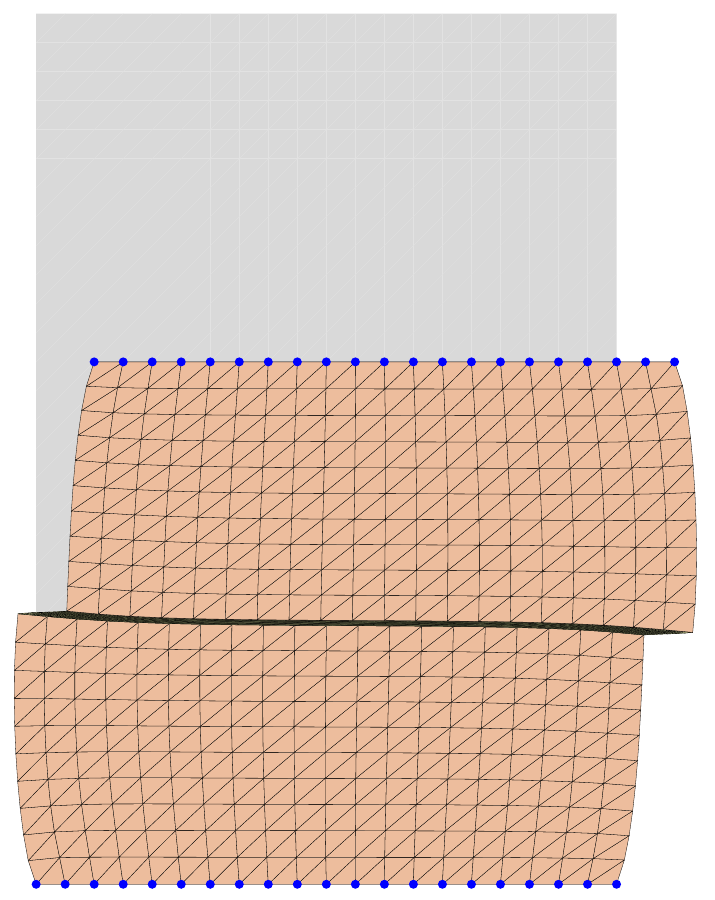}}
\put( 0.0,0.5){a) $\gamma=10^{-3}$}
\put( 4.,1.0){\includegraphics[width=0.24\textwidth]{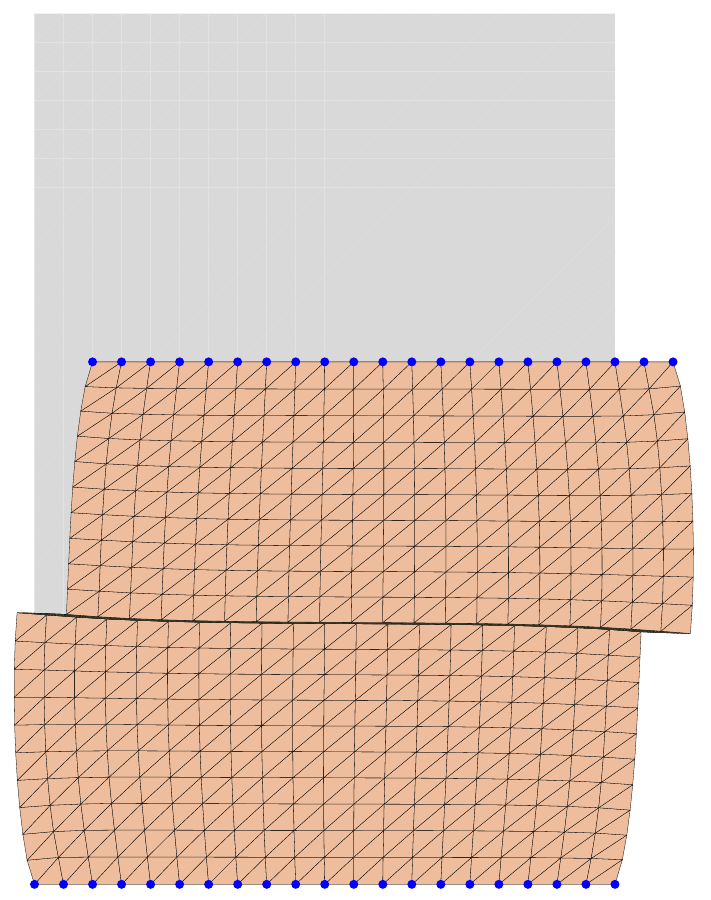}}
\put( 4.0,0.5){b) $\gamma=10^{-4}$}
\put( 8,1.0){\includegraphics[width=0.24\textwidth]{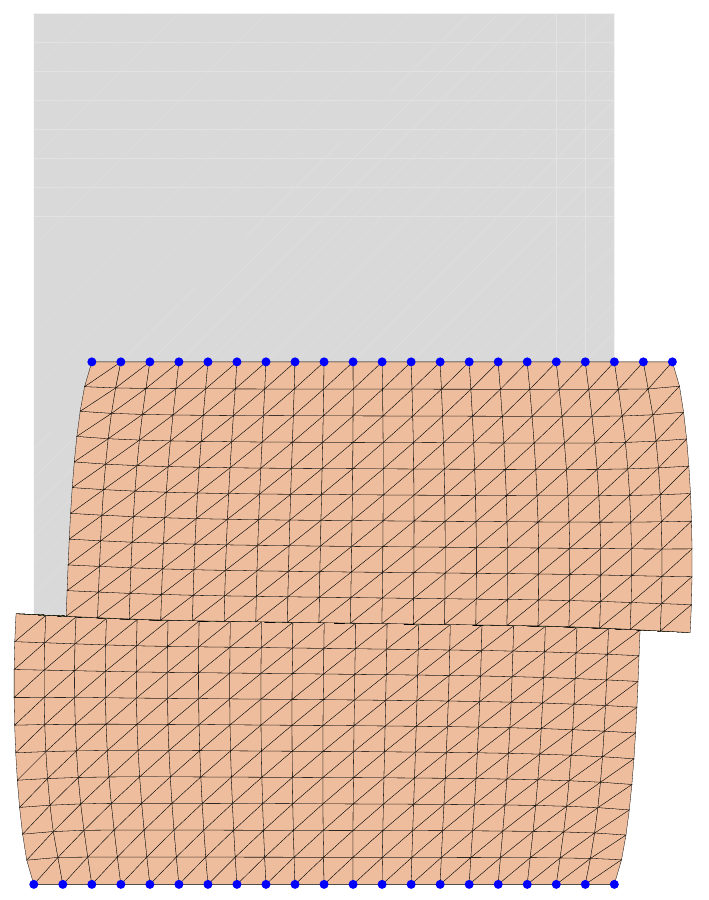}}
\put( 8.0,0.5){c) $\gamma=10^{-5}$}
\put(12,1.0){\includegraphics[width=0.24\textwidth]{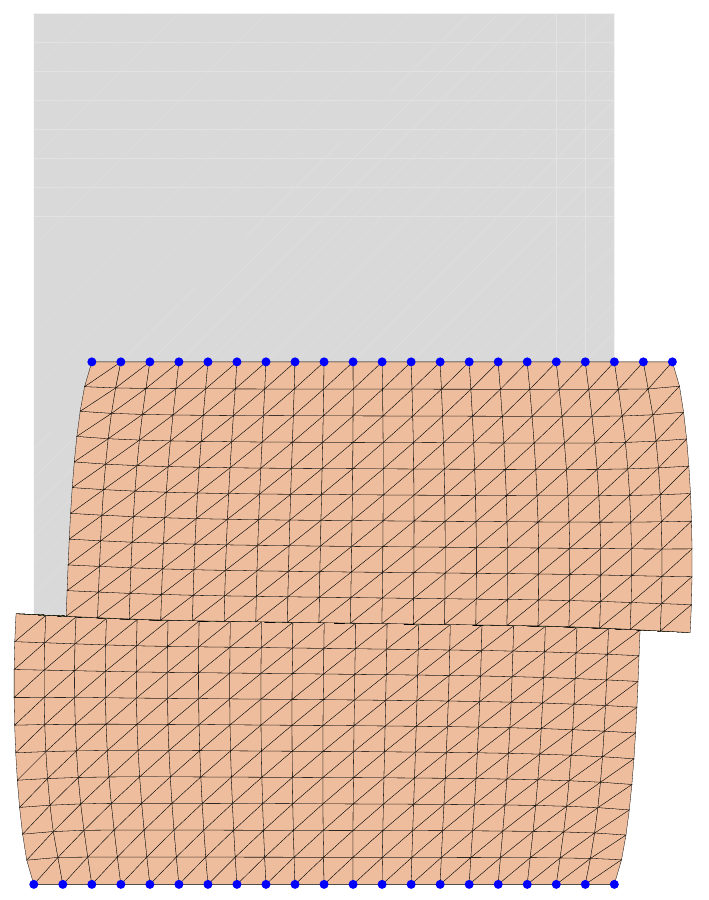}}
\put(12.0,0.5){d) $\gamma=10^{-8}$}
\end{picture}
\vspace*{-6mm}
\caption{
Final deformed configurations of the two-block benchmark under combined horizontal and vertical loading for different values of the third-medium scaling parameter $\gamma$. 
The gray overlay indicates the reference configuration.
}
\label{fig:contactBlocksVertDeformedShear}
\end{Figure}

The final deformed configurations for case~2 are shown in Fig.~\ref{fig:contactBlocksVertDeformedShear}. 
The superposed reference configuration illustrates the horizontal shift of the upper block and the resulting distortion of the third medium layer. 
Compared with the vertical loading case, the deformation state is visibly more asymmetric since the third medium has to accommodate both normal compression and tangential motion. 
For larger values of $\gamma$, the intermediate layer remains more pronounced and the two solid blocks stay further apart. 
For smaller values of $\gamma$, the interface closes more tightly while the deformation remains regular. 
The comparison confirms that the proposed stabilization controls the third-medium distortion also in the presence of shear-dominated contact kinematics.

\subsection{Self-contact within a box}
\textbf{Comparison of different interpolations.}
As a second benchmark, the proposed TMC formulation is tested with the classical self-contact-within-a-box problem shown in Fig.~\ref{fig:BoxT1T1}a) and b). 
A deformable solid frame encloses a third medium region inside the rectangular domain $[0,2]\times[0,0.5]$. 
Thus, the frame represents the physical solid, while the interior of the box is filled by the fictitious contact medium. 
For the solid material, the bulk and shear moduli are chosen as $K=20$ and $\mu=10$, respectively, following the benchmark setup in \cite{WriKorJun:2025:atm}. 
Boundary conditions are prescribed at three points: the lower left point is fixed by $\overline{\bu}(0,0)=\bzero$, the lower right point is constrained in vertical direction by $\overline{u}_2(2,0)=0$, and the upper loading point is displaced by $\overline{\bu}(1,0.5)=[0,-1]^T$. 
Inside the third medium, the auxiliary-field stabilization introduced above is applied. 
An unstructured triangular mesh serves as discretization. 
To assess the role of the interpolation order, the continuous pairs $\mathrm{T}_{1}^{u}\mathrm{T}_{1}^{\Theta}$, $\mathrm{T}_{2}^{u}\mathrm{T}_{2}^{\Theta}$ and $\mathrm{T}_{2}^{u}\mathrm{T}_{1}^{\Theta}$ are compared with the discontinuous variants $\mathrm{T}_{1}^{u}\mathrm{T}_{1}^{\Theta,d}$, $\mathrm{T}_{1}^{u}\mathrm{T}_{0}^{\Theta,d}$, and $\mathrm{T}_{2}^{u}\mathrm{T}_{1}^{\Theta,d}$. 
Special attention is given to the sensitivity of the contact response with respect to the penalty parameter $p_\Theta$ and the regularization parameter $\alpha_r$. 
As accuracy measure, the residual gap between the upper and lower parts of the deforming frame is evaluated, where smaller values indicate a tighter enforcement of the third-medium contact response.
The reference and deformed configurations of the self-contact-within-a-box benchmark are compared in Fig.~\ref{fig:BoxT1T1} for two penalty parameters. 
A coarse unstructured mesh with 354 elements is used in Figs.~\ref{fig:BoxT1T1}a) and c), whereas Fig.~\ref{fig:BoxT1T1}b) and d) shows a refined discretization with 3208 elements. 
Both parameter sets are chosen such that a stable contact response with a small residual gap is obtained. 
As expected for penalty-type couplings, the suitable magnitude of $p_\Theta$ is problem-dependent and has to be balanced against the third-medium stiffness, the mesh size and the interpolation order. 
The configurations therefore serve as a first qualitative assessment of the proposed stabilization. 
A systematic parameter study is carried out in the following to identify robust and accurate choices of $p_\Theta$ and $\alpha_r$ for the different interpolation pairs, with the aim of keeping the residual gap between the upper and lower flanges small.
The sensitivity study starts with the penalty parameter $p_\Theta$, as summarized in Tab.~\ref{tab:sensitivity_study_ptheta_cont}. 
For the continuous interpolation pairs, all simulations converge for $p_\Theta=10^{-1}$ over the considered range of $\gamma$. 
This value is therefore selected for the subsequent study of the gradient penalty $\alpha_r$, reported in Tab.~\ref{tab:sensitivity_study_alpha_r_cont}. 
For this parameter choice, all tested combinations converge, which indicates a robust coupling between the auxiliary field and the deformation gradient.

\begin{Figure}[H]
\begin{picture}(0,3.0)
\unitlength1cm
\put( 0.0,0.5){\includegraphics[width=8cm]{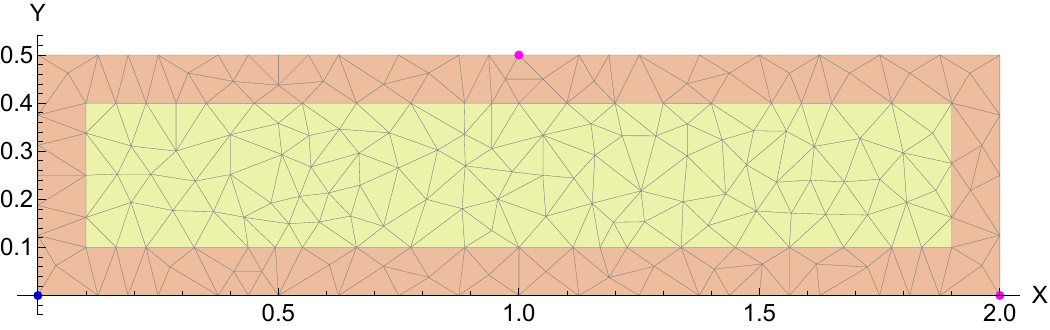}}
\put( 0.0,0.0){a)}
\put( 8.0,0.5){\includegraphics[width=8cm]{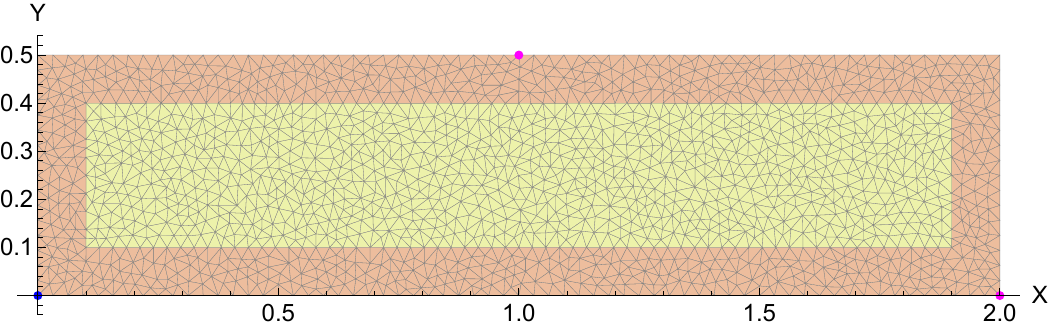}}
\put( 8.0,0.0){b)}
\put( 0.0,-4.8){\includegraphics[width=7.8cm]{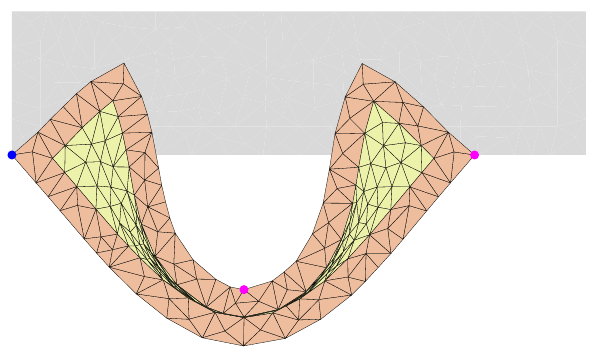}}
\put( 0.0,-3.5){c)}
\put( 8.0,-4.8){\includegraphics[width=7.8cm]{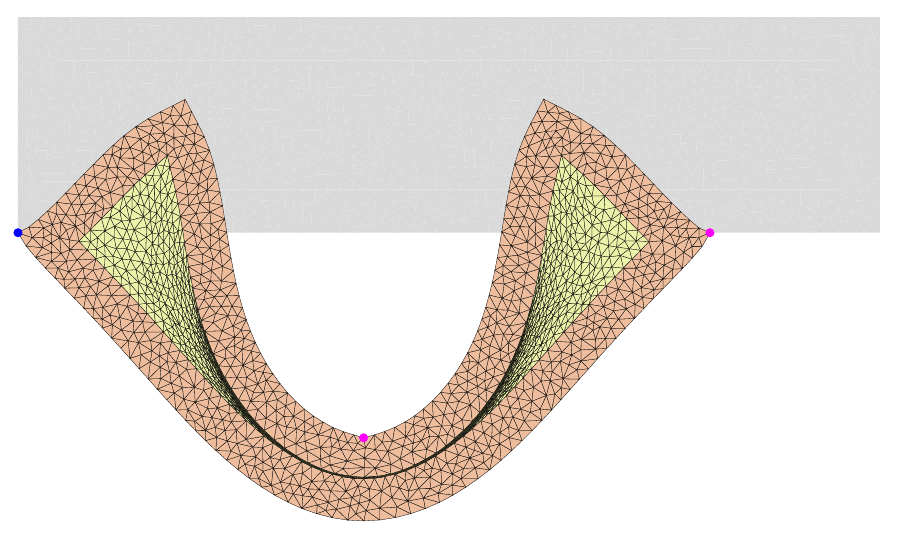}}
\put( 8.0,-3.5){d)}
\end{picture}
\vspace*{45mm}
\caption{
Reference and deformed configurations of the self-contact-within-a-box benchmark for two penalty parameters and the $\mathrm{T}_{1}^{u}\mathrm{T}_{1}^{\Theta}$ approximation. 
The solid frame is shown in orange, while the enclosed third medium is shown in yellow. 
Graphics a) and c) use $\gamma=10^{-5}$, $\alpha_r=10^{2}$ and $p_\Theta=10^{-2}$; graphics b) and d) use $\gamma=10^{-5}$, $\alpha_r=10^{2}$ and $p_\Theta=10^{-1}$. 
Panels a) and b) show the reference configurations, whereas c) and d) show the corresponding deformed configurations.
}
\label{fig:BoxT1T1}
\end{Figure}

The parameter studies in Tabs.~\ref{tab:sensitivity_study_ptheta_cont}, \ref{tab:sensitivity_study_alpha_r_cont} and~\ref{tab:sensitivity_study_ptheta_disc} show several consistent trends. 
For the continuous auxiliary-field interpolations, $p_\Theta=10^{-1}$ provides a robust parameter range for the present benchmark, since all tested values of $\gamma$ reach the final load step and produce small residual gaps. 
Larger penalty parameters do not necessarily improve the contact response. 
Instead, they can deteriorate the nonlinear convergence behavior and may lead to premature termination of the load path. 
The regularization parameter $\alpha_r$ mainly affects the residual gap up to a moderate value. 
Beyond this range, the additional improvement becomes small, indicating a saturation of the gradient-regularization effect for this benchmark. 
Comparing the interpolation pairs, the $\mathrm{T}_{2}^{u}\mathrm{T}_{1}^{\Theta}$ approximation gives results that are very close to those obtained with $\mathrm{T}_{2}^{u}\mathrm{T}_{2}^{\Theta}$, while using fewer auxiliary degrees of freedom. 
Thus, enriching the displacement field while keeping a lower-order auxiliary field appears to be an efficient choice for the considered self-contact problem. 
The $\mathrm{T}_{1}^{u}\mathrm{T}_{1}^{\Theta}$ approximation remains attractive from a low-order perspective, but its residual gaps are larger than those of the quadratic displacement formulations. 
The discontinuous auxiliary-field interpolations in Tab.~\ref{tab:sensitivity_study_ptheta_disc} show a different behavior. 
For the quadratic displacement approximation, the $\mathrm{T}_{2}^{u}\mathrm{T}_{1}^{\Theta,d}$ and $\mathrm{T}_{2}^{u}\mathrm{T}_{0}^{\Theta,d}$ variants give almost identical results over the investigated parameter range. 
This indicates that the additional local auxiliary degrees of freedom of the discontinuous linear field do not translate into a stronger regularization effect for this benchmark. 
Both variants provide small residual gaps for moderate values of $\gamma$, but the convergence deteriorates when the third medium becomes very soft. 
For example, at $p_\Theta=10^{-1}$ both quadratic discontinuous variants reach the final load step down to $\gamma=10^{-5}$, whereas smaller values of $\gamma$ lead to premature termination of the load path. 
The $\mathrm{T}_{1}^{u}\mathrm{T}_{1}^{\Theta,d}$ interpolation is less favorable in this example, since several parameter sets stop well before the final load step and the results are nearly insensitive to $p_\Theta$. 
The element-wise constant $\mathrm{T}_{1}^{u}\mathrm{T}_{0}^{\Theta,d}$ variant is more robust with respect to load completion for moderate penalty values, but it has to be interpreted as a limiting case because $\nabla\BTheta$ vanishes within each element. 
Consequently, it does not provide a true gradient-driven stabilization and mainly reflects the effect of the local penalty coupling. 
Overall, the behavior of the discontinuous variants is consistent with the one-dimensional discussion in Sec.~\ref{sec:1DRegularization}: element-local auxiliary fields are attractive from an implementation point of view, but the missing continuity-driven inter-element coupling weakens the gradient-type regularization mechanism compared with the continuous auxiliary-field formulations.

\textbf{Comparison with deformation-gradient averaging - simple box.}
A comparison with established regularization strategies is essential to assess both the robustness and the computational efficiency of the proposed formulation. 
Recently, Faltus~et~al.~\cite{FalAmaHor:2026:dga} introduced a deformation-gradient averaging approach that avoids additional degrees of freedom and therefore does not enlarge the global system of equations. 
In their formulation, a representative deformation gradient $\bar{\bF}$ is computed at the element center and enforced at the remaining integration points by a penalty contribution. 
In this way, variations of $\bF$ inside an element are suppressed without explicitly regularizing $\nabla\bF$. 
The method can therefore be interpreted as an element-wise deformation-gradient regularization.
For a fair comparison, the same benchmark and the same $Q_1$ displacement interpolation are used in both formulations. 
The formulation by Faltus~et~al. is evaluated with a $Q_1^u$ discretization, while the present mixed formulation uses a $Q_1^uQ_1^\Theta$ interpolation. 
Although the proposed approach introduces additional auxiliary degrees of freedom, it also provides a direct mixed representation of the deformation-gradient-like stabilization field. 
The parameters are kept fixed for all mesh refinements. 
For the proposed formulation, $p_\Theta=10$, $\alpha_r=1000$ and $\gamma=10^{-6}$ are used. 
For the deformation-gradient averaging approach, the penalty parameter is chosen as $\kappa_{\bar F}=0.5$ with $\gamma_{\bar F}=10^{-6}$.

\begin{table}[H]
\centering
\scriptsize
\setlength{\tabcolsep}{4pt}
\caption{
Performance comparison of the proposed mixed TMC formulation and the deformation-gradient averaging approach by Faltus et al.~\cite{FalAmaHor:2026:dga}. 
The proposed formulation uses a $Q_1^uQ_1^\Theta$ interpolation, whereas the formulation by Faltus et al. uses a $Q_1^u$ displacement interpolation. 
The reported total time is computed as the sum of residual-and-stiffness assembly time and solver time.
}
\begin{tabular}{c|r|r|r|r|r|r|r}
\hline\hline
formulation 
& elements 
& unknowns 
& load steps 
& iterations 
& K\&R time 
& solver time 
& total time \\
\hline\hline

\multirow{7}{*}{%
\begin{tabular}{c}
$Q_1^uQ_1^\Theta$\\[2mm]
proposed
\end{tabular}
}
& 100    & 552     & 37 & 244 & 0.159 s  & 0.178 s   & 0.337 s \\
& 400    & 1934    & 35 & 230 & 0.362 s  & 0.419 s   & 0.781 s \\
& 1600   & 7194    & 34 & 218 & 1.169 s  & 1.228 s   & 2.397 s \\
& 6400   & 27698   & 29 & 181 & 2.956 s  & 4.124 s   & 7.080 s \\
& 25600  & 108642  & 26 & 154 & 9.066 s  & 18.473 s  & 27.539 s \\
& 102400 & 430274  & 25 & 153 & 36.946 s & 107.640 s & 144.586 s \\
& 409600 & 1712514 & 30 & 353 & 315.820 s & 1579.500 s & 1895.320 s \\

\hline

\multirow{7}{*}{%
\begin{tabular}{c}
$Q_1^u$\\[2mm]
Faltus et al. \cite{FalAmaHor:2026:dga}
\end{tabular}
}
& 100    & 248     & 66  & 464 & 0.184 s  & 0.145 s   & 0.329 s \\
& 400    & 898     & 62  & 422 & 0.229 s  & 0.380 s   & 0.609 s \\
& 1600   & 3398    & 61  & 392 & 0.723 s  & 0.905 s   & 1.628 s \\
& 6400   & 13198   & 65  & 446 & 2.494 s  & 2.818 s   & 5.312 s \\
& 25600  & 51998   & 79  & 594 & 12.754 s & 19.501 s  & 32.255 s \\
& 102400 & 206398  & 113 & 868 & 62.589 s & 120.519 s & 183.108 s \\
& 409600 & 822398  & 108 & 893 & 234.361 s & 569.005 s & 803.366 s \\

\hline\hline
\end{tabular}
\label{tab:performance_faltus_vorwerk}
\end{table}

The results in Tab.~\ref{tab:performance_faltus_vorwerk} show the expected increase in system size caused by the additional auxiliary field of the proposed mixed formulation. 
For all mesh refinements, however, this larger system size is accompanied by a substantially more stable nonlinear load path. 
The proposed formulation consistently requires fewer load steps than the deformation-gradient averaging approach and also reduces the total number of Newton iterations for all considered meshes up to $409600$ elements. 
At $102400$ elements, the computation reaches the final configuration in $25$ load steps and $153$ iterations, compared with $113$ load steps and $868$ iterations for the formulation by Faltus~et~al.~\cite{FalAmaHor:2026:dga}. 
At this mesh level, the improved nonlinear robustness compensates for the additional algebraic cost and reduces the total time from $183.108\,\mathrm{s}$ to $144.586\,\mathrm{s}$. 
The largest discretization with $409600$ elements represents an extremely fine mesh and leads to more than $1.7$ million unknowns for the mixed formulation. 
At this scale, the cost of the additional auxiliary degrees of freedom becomes dominant and the total solution time increases to $1895.320\,\mathrm{s}$, compared with $803.366\,\mathrm{s}$ for the deformation-gradient averaging approach. 
Even in this demanding case, however, the proposed formulation reaches the final load in only $30$ load steps, whereas the reference approach requires $108$ load steps. 
Thus, the mixed auxiliary-field formulation trades a larger algebraic system for a considerably more robust load stepping behavior, with the overall efficiency depending on the balance between nonlinear robustness and the cost of the additional auxiliary unknowns.

\subsection{Self-contact of a C-shaped box}
The proposed TMC stabilization is further assessed with the classical C-shaped box benchmark shown in Fig.~\ref{fig:BVP_CShape}. 
A two-dimensional C-shaped hyperelastic solid with outer dimensions $1000\times500$ and a wall thickness of $100$ is considered. 
Inside the gap, a third medium fills the open space, while an additional third-medium layer of thickness $12.5$ is placed at the right boundary to allow contact between the upper and lower arms. 
Along the left boundary, the solid is fully clamped. 
At the upper right corner $P(1000,500)$, a prescribed displacement $\bar{\bu}=(0,\bar{u}_y)^T$ drives the closing motion of the C-shape. 
Large rotations, severe compression of the third medium and progressive self-contact are therefore induced by a simple displacement-controlled loading. 
For the solid, the compressible Neo-Hookean material in Eq.~\ref{eq:NeoHook} is used with $K_{\mathrm{s}}=5/3$ and $\mu_{\mathrm{s}}=5/14$.
As a first check, the influence of the discretization is assessed with the continuous $\mathrm{T}_{1}^{u}\mathrm{T}_{1}^{\Theta}$ interpolation. 
Fig.~\ref{fig:CShape_meshes} shows the deformation evolution for four increasingly refined meshes. 

\begin{Figure}[h]
\begin{picture}(0,6.8)
\unitlength1cm
\put( 0.0,1.0){\includegraphics[width=0.7\textwidth]{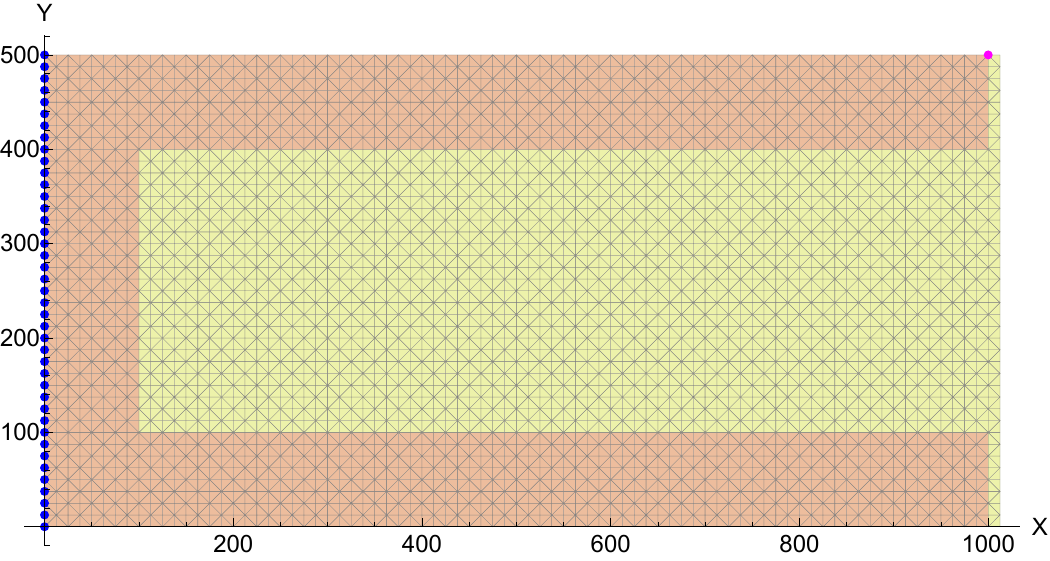}}
\put(11.0,6.0){\textbf{Parameters}}
\put(11.3,5.3){Shear modulus: $\mu=5/14$}
\put(11.3,4.8){Bulk modulus: $K=5/3$}
\end{picture}
\vspace*{-8mm}
\caption{
Reference configuration of the C-shaped box benchmark for finite deformation contact. 
Orange indicates the physical solid domain, and yellow indicates the third medium. 
Blue dots mark the clamped boundary, while the violet dot marks the point where the prescribed displacement is applied.
} 
\label{fig:BVP_CShape}
\end{Figure}

Since the penalty coupling is mesh-dependent, the parameter $p_\Theta$ is adjusted for each discretization rather than kept fixed. 
This procedure should therefore not be interpreted as a classical mesh-convergence study with one fixed parameter set, but as a calibrated mesh sequence used to assess whether the same deformation mechanism can be obtained over a range of discretizations. 
For the four meshes shown in Fig.~\ref{fig:CShape_meshes}, the values $p_\Theta=\{0.001,0.005,0.01,0.05\}$ are used. 
Across all meshes, the same global deformation mode is obtained, i.e., the upper arm bends into the cavity, the third medium is squeezed into a narrow band, and self-contact develops along the inner side of the C-shape. 
With increasing mesh resolution, the deformation path becomes smoother and the compressed third-medium layer is represented more sharply. 
Qualitative agreement across all calibrated meshes indicates that the proposed stabilization is not tied to a particular mesh density.

\begin{Figure}[h]
\begin{picture}(0,5.5)
\unitlength1cm
\put( 0.0,0.0){a) 323 DOF}
\put( 0.0,0.4){\includegraphics[width=4cm]{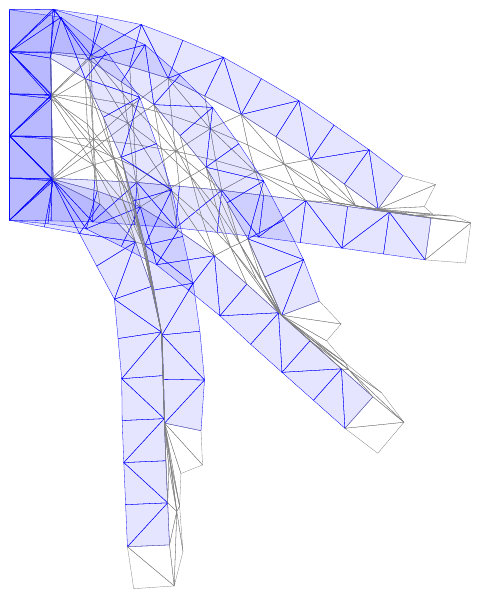}}
\put( 4.0,0.0){b) 1053 DOF}
\put( 4.0,0.4){\includegraphics[width=4cm]{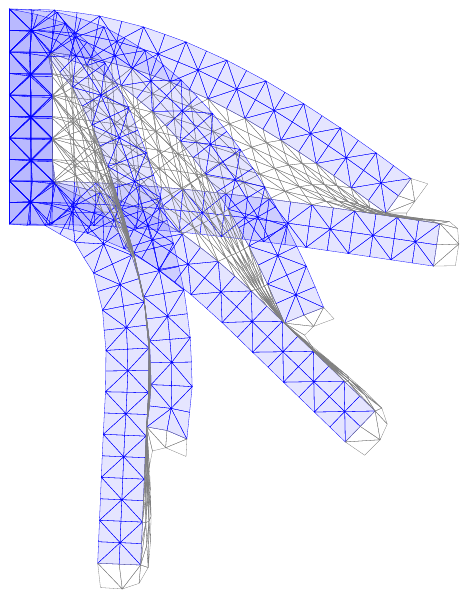}}
\put( 8.0,0.0){c) 3761 DOF}
\put( 8.0,0.4){\includegraphics[width=4cm]{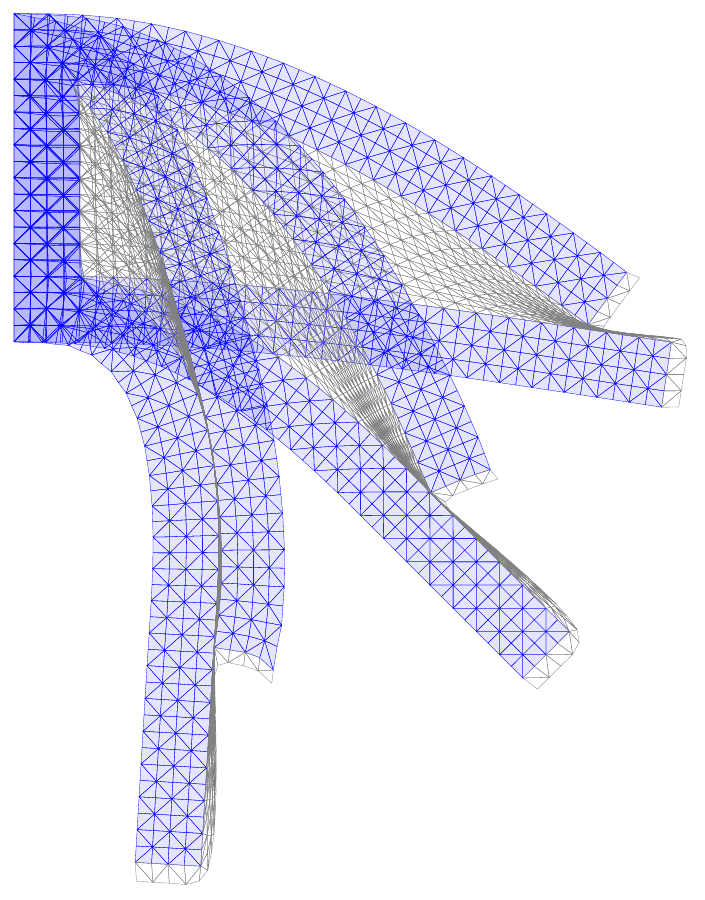}}
\put(12.0,0.0){d) 14169 DOF}
\put(12.0,0.4){\includegraphics[width=4cm]{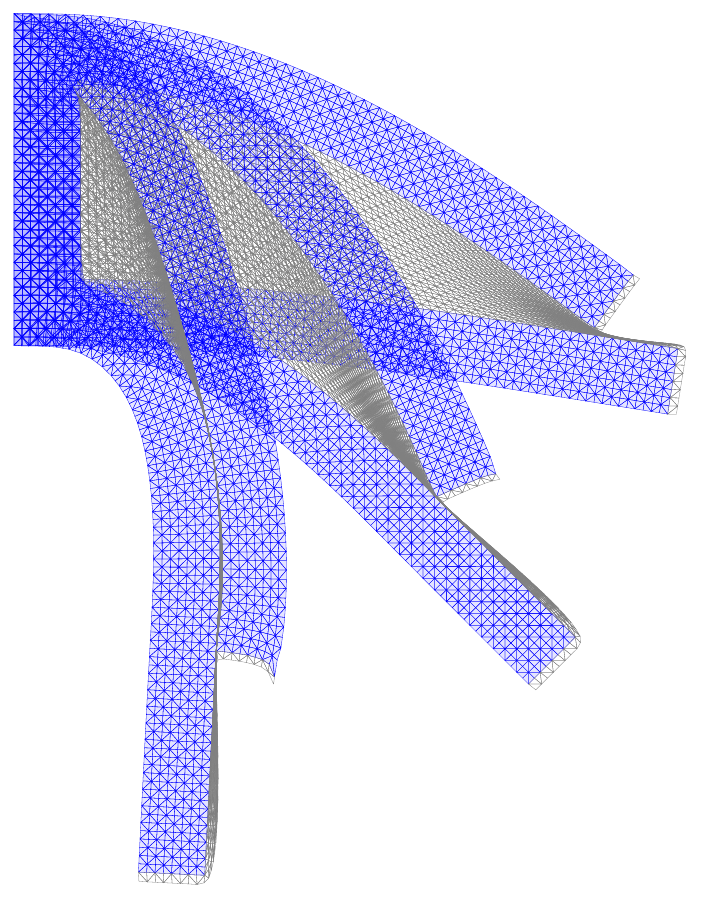}}
\end{picture}
\vspace*{1mm}
\caption{
Deformation history of the C-shaped benchmark for a calibrated sequence of increasingly refined meshes using the continuous $\mathrm{T}_{1}^{u}\mathrm{T}_{1}^{\Theta}$ discretization. 
The penalty parameter $p_\Theta$ is adjusted for each mesh. 
Snapshots correspond to the load parameters $\lambda_{\mathrm{load}}=0.4$, $\lambda_{\mathrm{load}}=0.7$ and $\lambda_{\mathrm{load}}=1.0$.
}
\label{fig:CShape_meshes}
\end{Figure}

\begin{Figure}[h]
\begin{picture}(0,6.8)
\unitlength1cm
%\put( 0.0,7.2){T$_1^u$T$_1^\Theta$ discretization}
\put( 0.0,4.1){\includegraphics[width=0.33\textwidth]{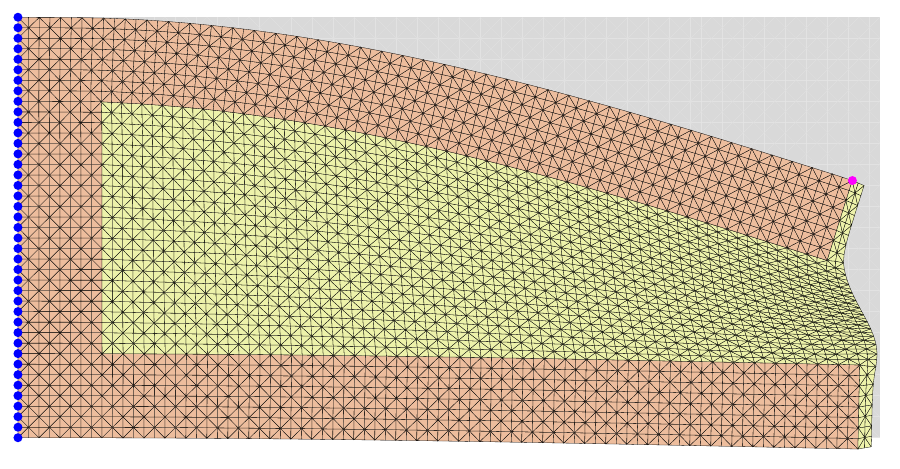}}
\put( 0.0,3.9){a) $\lambda_{\mathrm{load}}=0.2$}
\put( 0.0,0.0){\includegraphics[width=0.33\textwidth]{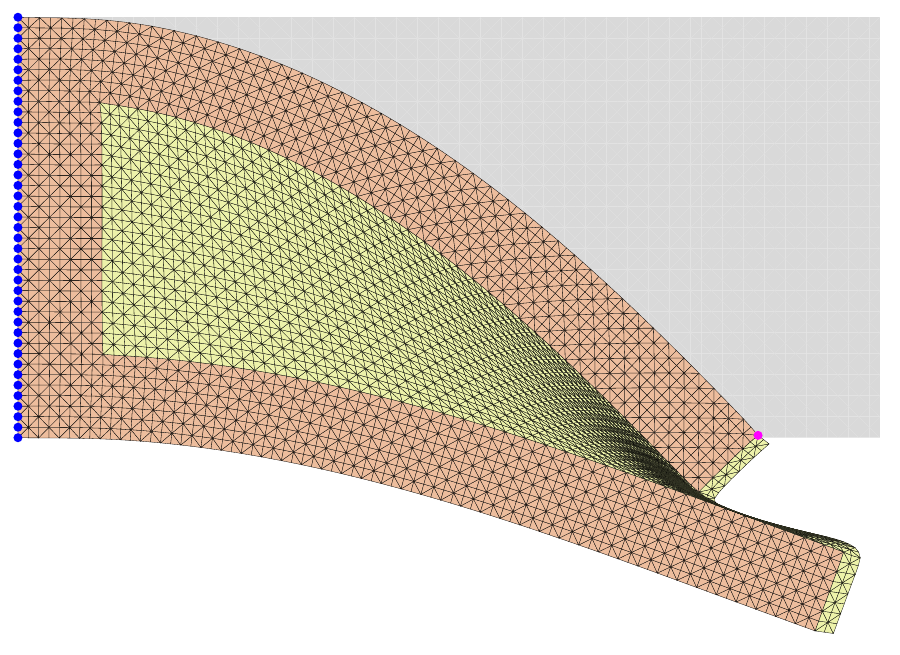}}
\put( 0.0,0.5){b) $\lambda_{\mathrm{load}}=0.5$}
\put( 5.5,1.7){\includegraphics[width=0.33\textwidth]{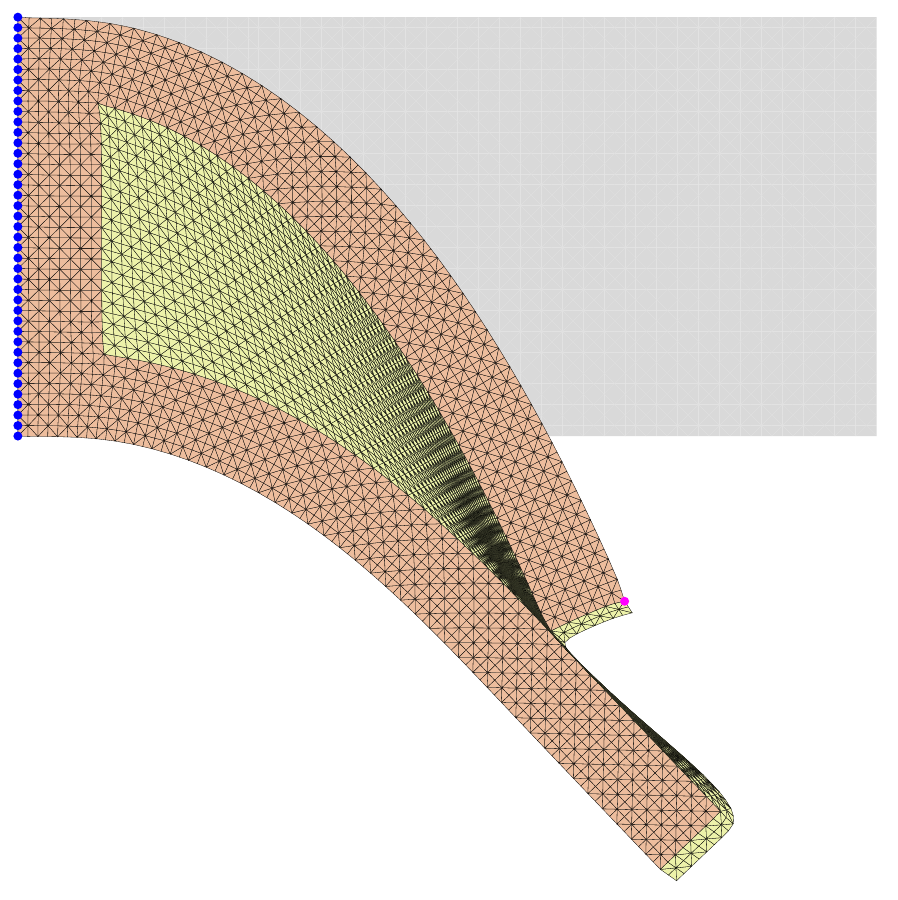}}
\put( 5.5,4.05){c)}
\put( 8.5,6.0){$\lambda_{\mathrm{load}}=0.7$}
\put(11,0.2){\includegraphics[width=0.33\textwidth]{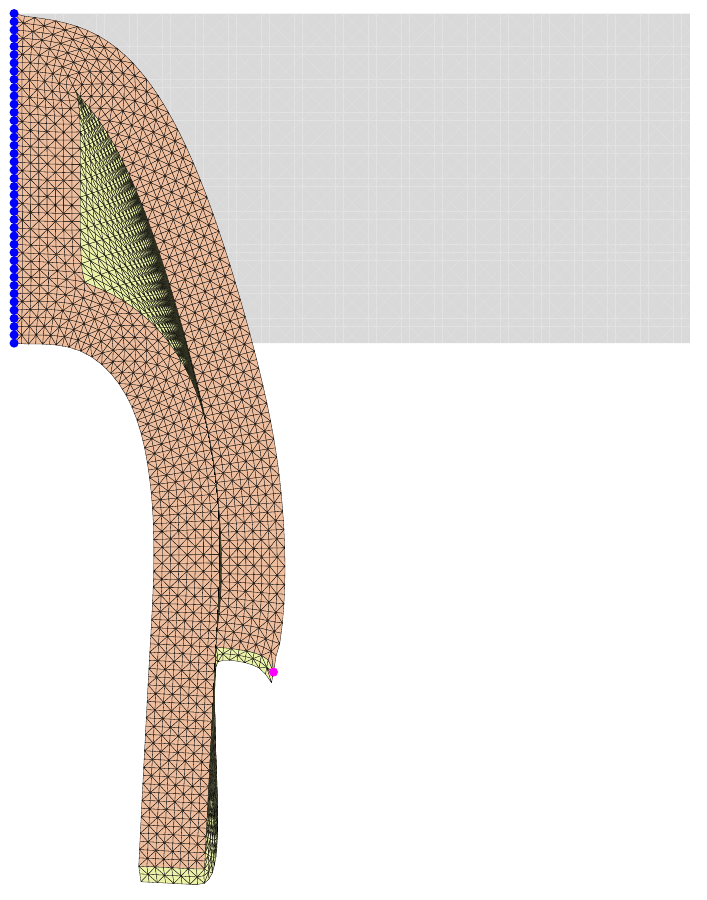}}
\put(11.0,4.0){d)}
\put(13.0,6.0){$\lambda_{\mathrm{load}}=1.0$}
\put( 6.0,2.0){\textbf{Parameters}}
\put( 6.0,1.5){$\gamma=10^{-8}$}
\put( 6.0,1.0){$\alpha_r=10^0$}
\put( 6.0,0.5){$p_\Theta=5\cdot 10^{-2}$}
\end{picture}
\vspace*{-2mm}
\caption{
Deformation history of the C-shaped benchmark using the continuous $\mathrm{T}_{1}^{u}\mathrm{T}_{1}^{\Theta}$ discretization. 
Snapshots are shown for increasing load parameters $\lambda_{\mathrm{load}}$. 
All deformations are plotted without displacement scaling. 
The gray overlay indicates the reference configuration.
}
\label{fig:CShapeT1T1}
\end{Figure}

A more detailed view of the deformation process is given in Fig.~\ref{fig:CShapeT1T1}. 
At early loading, the upper arm starts to rotate downward while the third medium remains broadly distributed inside the cavity. 
With increasing load, the gap closes and the third medium is progressively compressed between the approaching solid surfaces. 
At the final load step, deformation localizes along the contact region, while the third-medium mesh remains regular and no element collapse is observed. 
Even with the continuous $\mathrm{T}_{1}^{u}\mathrm{T}_{1}^{\Theta}$ discretization, the large rotation and subsequent self-contact are captured robustly and the overall deformation pattern remains smooth and stable throughout the loading path.

\begin{Figure}[h]
\begin{picture}(0,7.5)
\unitlength1cm
\put( 0.0,2.0){\rotatebox{90}{Vertical reaction force $R_y$}}
\put( 2.0,0.5){Prescribed displacement $u_y$}
\put( 0.5,0.5){a)}
\put( 0.5,1.0){\includegraphics[width=7.5cm]{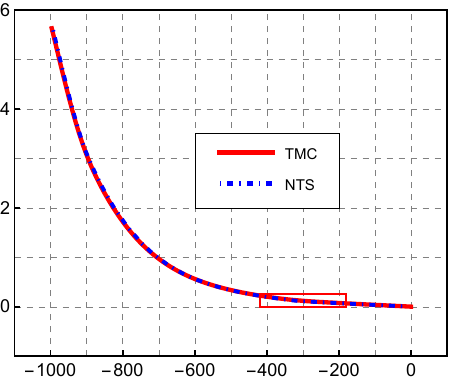}}
\put(10.0,0.5){Prescribed displacement $u_y$}
\put( 8.5,0.5){b)}
\put( 8.5,1.0){\includegraphics[width=7.5cm]{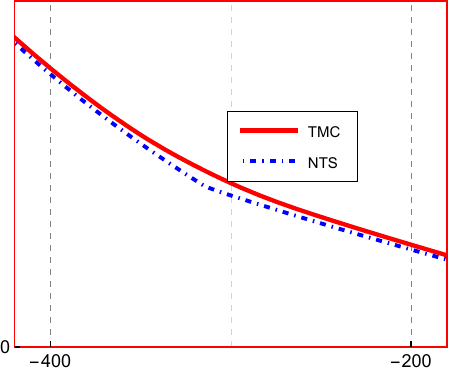}}
\end{picture}
\vspace*{-3mm}
\caption{
Reaction-force-displacement response of the C-shaped benchmark. 
The third medium contact formulation is compared with a classical node-to-segment contact formulation. 
a) Full loading path. 
b) Magnified view of the contact regime.
}
\label{fig:CShape_LoadDisp}
\end{Figure}

Fig.~\ref{fig:CShape_LoadDisp} shows the vertical reaction force as a function of the prescribed vertical displacement $u_y$. 
The reaction force is obtained from the vertical support reaction associated with the imposed displacement at the loading point. 
Over most of the loading range, the TMC curve follows the NTS reference closely. 
During the initial bending-dominated phase, both curves are almost indistinguishable and remain in good agreement after contact is established. 
Only in the magnified initial contact regime in Fig.~\ref{fig:CShape_LoadDisp}b, a small deviation becomes visible. 
Such a deviation is expected for a continuum contact regularization with finite third-medium stiffness. 
For engineering-scale simulations, the difference remains small compared with the overall force level, while the TMC formulation avoids explicit contact search and provides a smooth transition into self-contact.

\textbf{Comparison with deformation-gradient averaging - C-box.}
A second comparison with the deformation-gradient averaging approach of Faltus~et~al.~\cite{FalAmaHor:2026:dga} is performed for the C-shaped benchmark. 
Here, the focus is not on computational cost, but on the robustness of the third-medium deformation under different interpolation choices. 
Fig.~\ref{fig:CShape_Faltus} shows the deformed configurations at the load parameter $\lambda_{\mathrm{load}}=0.2$. 
$\bar \bF$ is calculated at the center points of each individual element.
A similar approach calculating the volume averaged $\bar \bF=\frac{1}{V}\int_{\B}\bF$dv reaches the same results.

\begin{Figure}[H]
\begin{picture}(0,0.2)
\unitlength1cm
\put( 0.0,0.1){\textbf{Continuous}}
\put( 5.5,0.1){\textbf{Discontinuous}}
\put(11.0,0.1){\textbf{Faltus et al.}}
\put( 0.0,-2.7){\includegraphics[width=5cm]{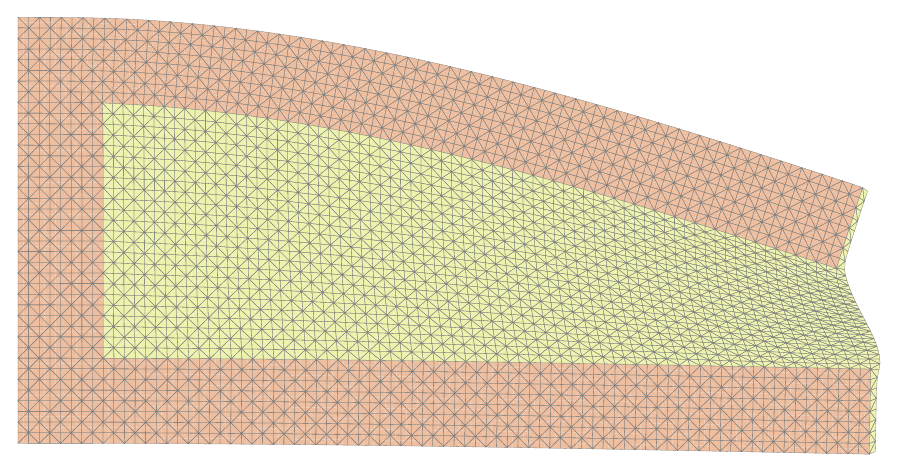}}
\put( 0.0,-3.0){a) $\mathrm{T}_{2}^{u}\mathrm{T}_{1}^{\Theta}:$}
\put( 0.0,-3.5){\phantom{a)}  $p_\Theta=10^{-2}$, $\alpha_r=1$}
\put( 5.5,-2.7){\includegraphics[width=5cm]{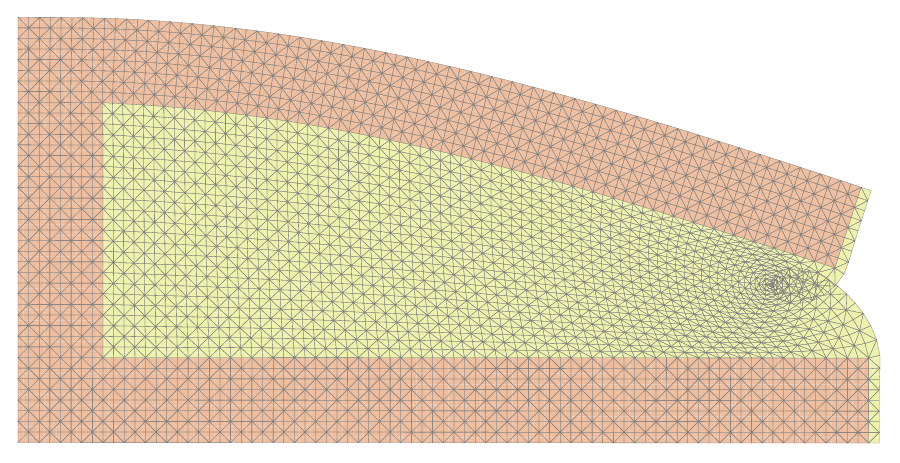}}
\put( 5.5,-3.1){b) $\mathrm{T}_{2}^{u}\mathrm{T}_{1}^{\Theta,d}:$ $p_\Theta=10^{-2}-10^{6}$}
\put(11.0,-2.7){\includegraphics[width=5.5cm]{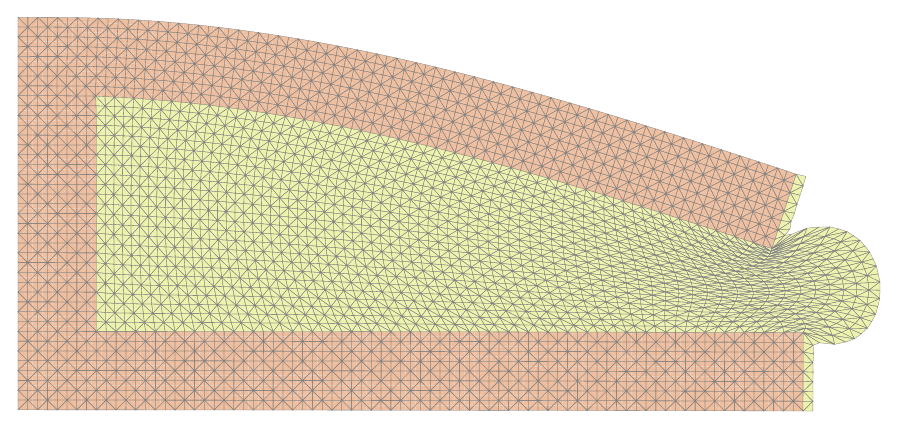}}
\put(11.0,-3.0){c) $\mathrm{T}_{2}^{u}:$ $\kappa=1$}
\put( 0.0,-6.7){\includegraphics[width=5cm]{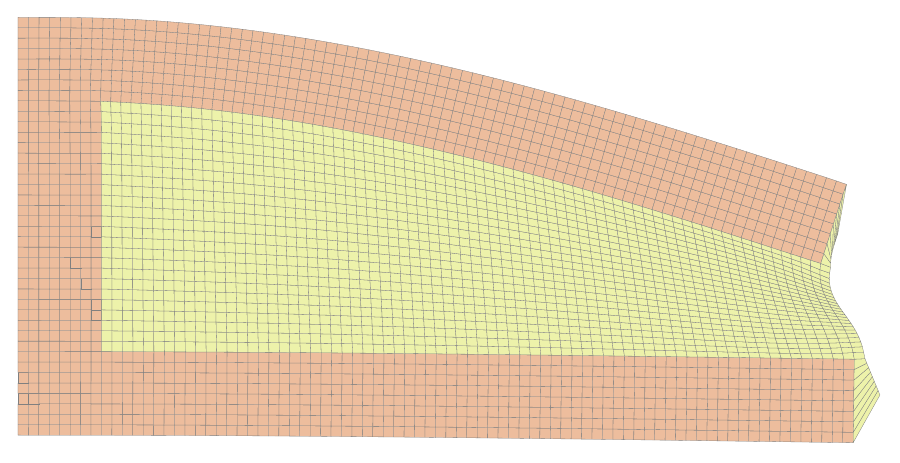}}
\put( 0.0,-7.0){d) $\mathrm{Q}_{1}^{u}\mathrm{Q}_{1}^{\Theta}:$}
\put( 0.0,-7.5){\phantom{b)} $p_\Theta=10^{-1}$, $\alpha_r=10^{-2}$}
\put( 5.5,-6.7){\includegraphics[width=5cm]{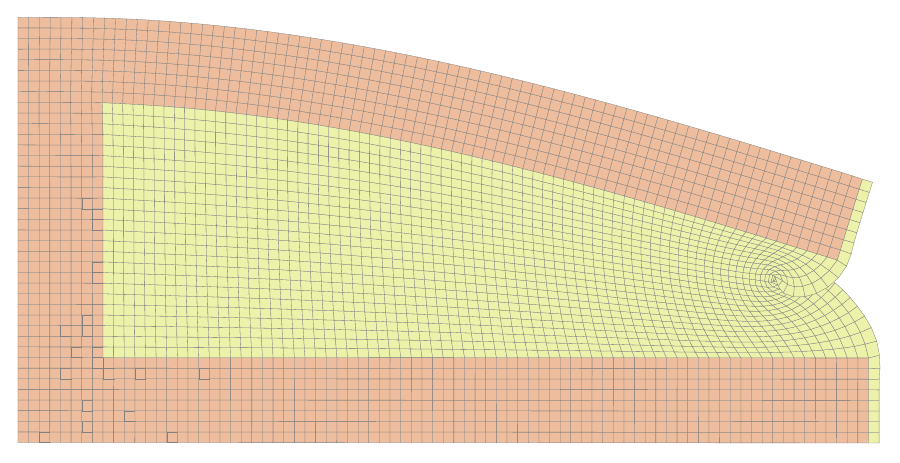}}
\put( 5.5,-7.1){e) $\mathrm{Q}_{1}^{u}\mathrm{Q}_{1}^{\Theta,d}:$ $p_\Theta=10^{-2}-10^{6}$}
\put(11.0,-6.7){\includegraphics[width=5cm]{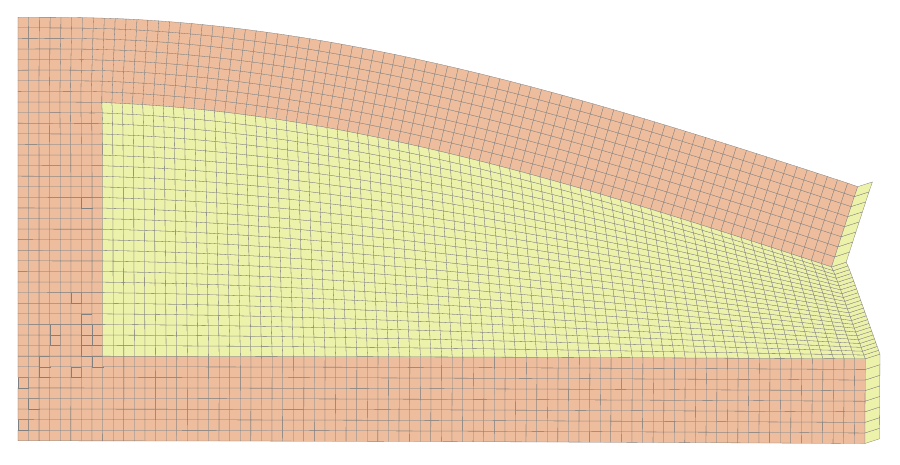}}
\put(11.1,-7.0){f) $\mathrm{Q}_{1}^{u}:$ $\kappa=10^{-2}$}
\end{picture}
\vspace*{75mm}
\caption{
Comparison of third-medium deformations in the C-shaped benchmark at $\lambda_{\mathrm{load}}=0.2$. 
The left column shows the proposed formulation with continuous auxiliary-field interpolations, the middle column shows the corresponding element-wise discontinuous variants, and the right column shows the deformation-gradient averaging approach of Faltus~et~al.~\cite{FalAmaHor:2026:dga}. 
The upper row shows triangular discretizations using a) $\mathrm{T}_{2}^{u}\mathrm{T}_{1}^{\Theta}$, b) $\mathrm{T}_{2}^{u}\mathrm{T}_{1}^{\Theta,d}$ and c) $\mathrm{T}_{2}^{u}$, while the lower row shows quadrilateral discretizations using d) $\mathrm{Q}_{1}^{u}\mathrm{Q}_{1}^{\Theta}$, e) $\mathrm{Q}_{1}^{u}\mathrm{Q}_{1}^{\Theta,d}$ and f) $\mathrm{Q}_{1}^{u}$. 
}
\label{fig:CShape_Faltus}
\end{Figure}

The continuous auxiliary-field formulation provides the most robust response in this comparison. 
For both triangular and quadrilateral discretizations, the third medium remains confined to the closing gap and follows the deformation of the C-shaped structure without visible element collapse. 
Although the result of the element $\mathrm{T}_{1}^{u}\mathrm{T}_{1}^{\Theta}$ is not displayed, the same behavior is observed in corresponding simulations.

The element-wise discontinuous variants behave fundamentally differently. 
Since $\nabla\BTheta$ vanishes inside each element for discontinuous low-order auxiliary interpolations, the gradient regularization does not provide an inter-element stabilization mechanism. For linear triangles the difference between $\boldsymbol F$ and ${\boldsymbol \Theta}$ is zero since the deformation gradient $\boldsymbol F$ is constant in the element. Thus the regularization cannot work. Furthermore, it is observed that for quadratic triangles the regularization also does not work. This is due to the fact that the displacement $\boldsymbol u_h$ is nearly a self-affine function for T$_2$ triangle. We note, that a self-affine (linear) displacement field occurs in the quadratic triangle when the displacements at the mid points of the element are given by
\begin{equation}
\bu_4=\frac{\bu_1+\bu_2}{2}\,,\quad
\bu_5=\frac{\bu_2+\bu_3}{2}\,\quad \mbox{and}\,\,\,
\bu_6=\frac{\bu_3+\bu_2}{2}.
\label{eq:mid_nodes}
\end{equation}
The resulting linear displacement field yields a constant deformation gradient $\boldsymbol F$ and thus, the direct stabilization term Eq.~\ref{eq:TMC_direct_gradF} is zero or  the difference in Eq.~\ref{eq:TMC_penalty_theta_F} is zero and the regularization does not work. 
Furthermore, it can be shown by a Taylor series expansion that, e.g.,  for the displacements of a  mid node $4$ of the element at  $\bx_m$ which lies between the vertex nodes $\bx_i$ and   $\bx_j$ that  the following estimate holds with $\bh=\frac12(\bx_j -\bx_i)$ for the quadratic displacement field
\begin{equation}
\begin{aligned}
\bu(\bx_m+\bh)
&=
\bu(\bx_m)
+
\nabla\bu(\bx_m)\,\bh
+
\frac{1}{2}
\nabla^2\bu(\bx_m)
[\bh,\bh]
+
\mathcal{O}(\|\bh\|^3)\\
\bu(\bx_m-\bh)
&=
\bu(\bx_m)
-
\nabla\bu(\bx_m)\,\bh
+
\frac{1}{2}
\nabla^2\bu(\bx_m)
[\bh,\bh]
+
\mathcal{O}(\|\bh\|^3)\,.
\end{aligned}
\end{equation}
Adding the two equations and solving for the mid point displacement yields
\begin{equation}
\bu(\bx_m)
=
\frac{\bu_i+\bu_j}{2}
-
\frac{1}{2}
\nabla^2\bu(\bx_m)
[\bh,\bh]
+
\mathcal{O}(\|\bh\|^4)
\end{equation}
and therefore, to second-order accuracy the displacement at the mid nodes is provided by one half of the sum of the vertex nodes. This means that Eq.~\ref{eq:mid_nodes} is fulfilled up to second order accuracy. We note that with $\mathbf h = \frac{h}{2}\mathbf t$  we obtain 
\begin{equation}
\frac12\nabla^2\bu(\mathbf x_m)[\bh,\bh]
=
\frac{h^2}{8}
\frac{\partial^2\bu}{\partial s^2}(\mathbf x_m).
\end{equation}
where $\mathbf t$ is the unit tangent vector at the edge and $s$ the coordinate along the edge.

For the quadratic triangular element, $\bu_m$ is still an independent nodal displacement. It is generally not exactly the average of the end-node displacements. But the difference between $\bu_m$ and the average is of order $h^2$, so it vanishes under mesh refinement. It is actually zero for an affine (linear) displacement field where $ \nabla^2  \bu_h=\mathbf 0$. Thus the effect of the regularisation functional $W^{\textrm{tm}}_p$ is considerably small and goes to zero for finer meshes and does not  work as regularization.
Even changing the penalty parameter $p_\Theta$ cannot prevent the observed element collapse and local self-penetration, and the parameter $\alpha_r$ has no effective influence in this setting. 
The deformation-gradient averaging yields the same result. 

Hence, for the triangular discretization, see Fig. \ref{fig:CShape_Faltus} c), the third medium is pushed laterally out of the gap and forms an artificial extrusion at the right boundary, whereas the quadrilateral case remains stable. 
Overall, the comparison shows that the stabilizing effect of the proposed method relies on a continuous auxiliary-field interpolation when low-order gradient information is needed. 
Discontinuous variants remove this mechanism, while deformation-gradient averaging remains more dependent on the underlying element type.

\subsection{Three-dimensional self-contact within a box}
Finally, the self-contact-within-a-box benchmark is extended to a three-dimensional setting in order to demonstrate that the proposed formulation is not restricted to plane problems. 
Following the two-dimensional setup, a deformable solid frame encloses a third medium region, and the imposed displacement drives the structure into self-contact. 
A low-order $\mathrm{T}_{1}^{u}\mathrm{T}_{0}^{\Theta,d}$ tetrahedral discretization with an element-wise auxiliary field is used in the third medium and selected material parameters $\gamma=10^{-5}$, $\alpha_r=1$ and $p_\Theta=1$.
Fig.~\ref{fig:BoxT1T0_3D} shows the deformed configuration from two perspectives. 
During loading, the upper part of the frame bends into the box, the third medium is compressed between the approaching solid surfaces, and the deformation remains regular. 
Overall, the example confirms that the auxiliary-field stabilized third medium formulation carries over directly to three-dimensional finite deformation contact problems.

\begin{Figure}[H]
\begin{picture}(0,4.0)
\unitlength1cm
\put( 0.0,0.0){\includegraphics[width=8cm]{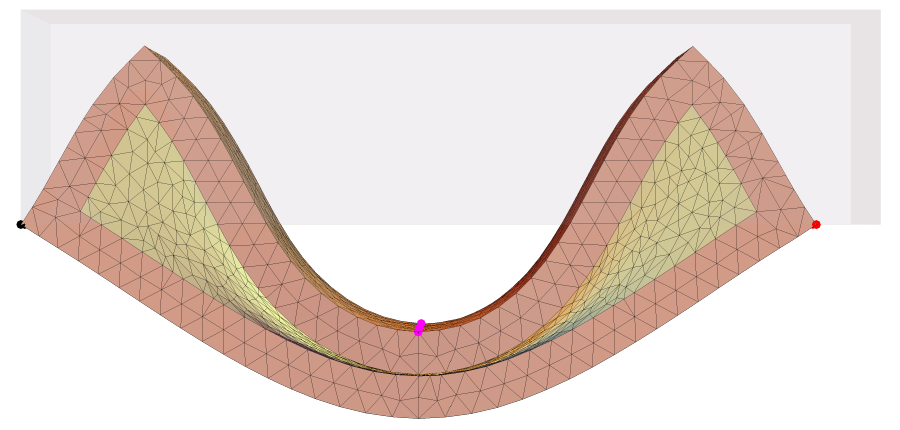}}
\put( 0.0,0.0){a)}
\put( 9.0,-1.5){\includegraphics[width=7cm]{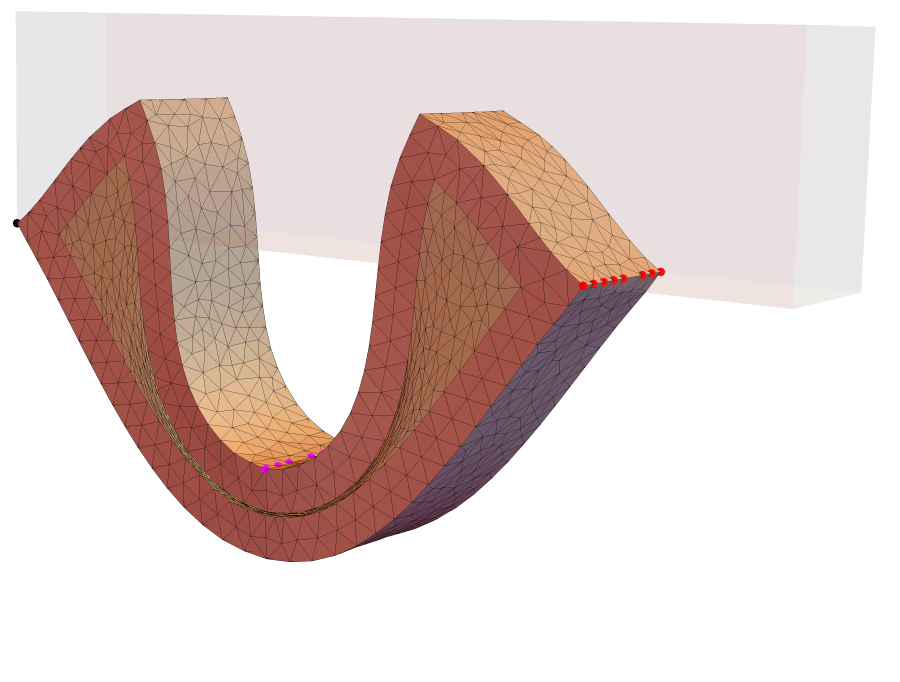}}
\put( 9.0,0.0){b)}
\end{picture}
\vspace*{5mm}
\caption{
Three-dimensional self-contact-within-a-box benchmark, with the physical solid frame enclosing the lighter third-medium region. 
The deformed configuration is shown from two perspectives. 
All deformations are plotted without displacement scaling.
}
\label{fig:BoxT1T0_3D}
\end{Figure}

\section{Conclusion and outlook}
In this work, an auxiliary-field stabilized mixed finite element formulation for third medium contact at finite deformations is presented.
The main idea is to introduce a deformation-gradient-like field $\BTheta$ in the third medium and to couple it weakly to the physical deformation gradient $\bF$ by a penalty contribution. 
Regularization acts on $\nabla\BTheta$ rather than on $\nabla\bF$ directly. 
In this way, gradient-type regularization is incorporated without evaluating second derivatives of the displacement field, which makes the approach applicable to first- and second-order finite elements.

A central observation is that the regularization effect depends strongly on the interpolation of the auxiliary field. 
For continuous auxiliary-field interpolations, neighboring elements share the nodal values of $\BTheta$. 
Element-wise changes of the deformation gradient are therefore transferred into spatial variations of the auxiliary field, which are penalized by the regularization parameter $\alpha_r$. 
Penalizing these spatial variations stabilizes the third medium and reduces the residual gap in the considered contact benchmarks. 
Element-wise discontinuous auxiliary fields are attractive from an implementation point of view, since their additional unknowns remain local to the element. 
However, the numerical results show that the missing continuity-driven inter-element coupling weakens the gradient-type stabilization mechanism compared with the continuous auxiliary-field formulations.

Across the numerical examples, large deformation contact, severe third-medium compression and progressive self-contact are captured in a stable manner.
In the two-block benchmark, the scaling parameter $\gamma$ controls the trade-off between residual gap and nonlinear solution effort, i.e., smaller values improve contact closure, while the Newton process becomes more demanding. 
The self-contact-within-a-box example confirms that moderate penalty parameters provide a robust coupling between $\BTheta$ and $\bF$, whereas overly large values of $p_\Theta$ may deteriorate convergence. 
For quadratic displacement approximations, the $\mathrm{T}_{2}^{u}\mathrm{T}_{1}^{\Theta}$ interpolation gives results close to the $\mathrm{T}_{2}^{u}\mathrm{T}_{2}^{\Theta}$ formulation while using fewer auxiliary degrees of freedom. 
The C-shaped benchmark further shows that the proposed third medium formulation gives a force-displacement response close to a classical augmented-Lagrange node-to-segment contact formulation, while retaining the main advantage of third medium contact, since no explicit contact search, active-set strategy or predefined contact interface is required.

Overall, the formulation combines the smooth continuum character of third medium contact with a low-order-compatible stabilization mechanism. 
Contact forces are transmitted through a fictitious continuum region, and self-contact develops naturally as part of the deformation process. 
Such a continuum-based treatment is particularly attractive for problems in which potential contact zones are not known in advance, for example in strongly deforming structures, contact-aided mechanisms and metamaterials with internal self-contact.

Several aspects remain open for future work. 
A systematic parameter selection strategy for $\gamma$, $p_\Theta$ and $\alpha_r$ is an important next step, especially for complex three-dimensional applications. 
Further work should address adaptive choices of the third-medium stiffness and more efficient solution strategies for large-scale problems. 
From an application point of view, the method appears especially promising for the analysis and design of mechanical metamaterials that undergo large shape changes and repeated internal contact.

\section*{Acknowledgement}
The authors gratefully acknowledge the computing time granted by the Center for Computational Sciences and Simulations (CCSS) of the University of Duisburg-Essen and provided on the supercomputer magnitUDE at the Zentrum f\"ur Informations- und Mediendienste (ZIM).

% === list of references
% ------- layout-datei --------------
%\bibliographystyle{plainnat}
\bibliographystyle{unsrtnat}
%\bibliographystyle{plaindin}
%\bibliographystyle{plaindin_shortname2}
%\bibliographystyle{elsarticle-num-names}
% ------- bib-datei --------------
\let\bibfont\relax
\bibliography{Contact_04}

%=== appendix
\begin{appendix}
\section{Appendix}

\begin{Table}[H]
\centering
\scriptsize
\setlength{\tabcolsep}{3pt}
\caption{
Residual gap, total Newton iteration count and achieved load parameter for the self-contact-within-a-box benchmark using continuous auxiliary-field interpolations.
Results are listed for different combinations of the third-medium scaling parameter $\gamma$ and the penalty parameter $p_\Theta$, while $\alpha_r=1$.
Red iteration counts denote diverged parameter sets.
}
\begin{tabular}{c|c|c|c||c|c|c||c|c|c}
  \hline\hline
  \multicolumn{4}{c||}{T$_2^u$T$_2^\Theta$} &
  \multicolumn{3}{c||}{T$_2^u$T$_1^\Theta$} &
  \multicolumn{3}{c}{T$_1^u$T$_1^\Theta$} \\[1mm]
  \hline

  \multicolumn{4}{l||}{$p_\Theta=10^{-2}$} &
  \multicolumn{3}{c||}{} &
  \multicolumn{3}{c}{} \\[1mm]
  \hline
  $\gamma$ & gap & iter. & $\lambda$ &
  gap & iter. & $\lambda$ &
  gap & iter. & $\lambda$ \\[1mm]
  \hline
  $10^{-3}$ & 0.0294122 & 118 & 1.000 & 0.0294120 & 118 & 1.000 & 0.0297254 & 118 & 1.000 \\
  $10^{-4}$ & 0.00790523 & 123 & 1.000 & 0.00790517 & 121 & 1.000 & 0.00877519 & 126 & 1.000 \\
  $10^{-5}$ & 0.00196967 & 128 & 1.000 & 0.00196965 & 126 & 1.000 & 0.00240124 & 164 & 1.000 \\
  $10^{-6}$ & 0.00051168 & 174 & 1.000 & 0.000511232 & 172 & 1.000 & 0.000673611 & 428 & 1.000 \\
  $10^{-7}$ & 0.000174526 & 388 & 1.000 & 0.000174144 & 409 & 1.000 & 0.000599197 & \textcolor{red}{274} & 0.422 \\
  $10^{-8}$ & 0.000093771 & 1098 & 1.000 & 0.0000937888 & 1025 & 1.000 & 0.000343603 & \textcolor{red}{337} & 0.436 \\
  \hline

  \multicolumn{4}{l||}{$p_\Theta=10^{-1}$} &
  \multicolumn{3}{c||}{} &
  \multicolumn{3}{c}{} \\[1mm]
  \hline
  $\gamma$ & gap & iter. & $\lambda$ &
  gap & iter. & $\lambda$ &
  gap & iter. & $\lambda$ \\[1mm]
  \hline
  $10^{-3}$ & 0.0251904 & 119 & 1.000 & 0.0251813 & 119 & 1.000 & 0.0297878 & 118 & 1.000 \\
  $10^{-4}$ & 0.00600028 & 128 & 1.000 & 0.00599759 & 128 & 1.000 & 0.0091427 & 125 & 1.000 \\
  $10^{-5}$ & 0.00141917 & 155 & 1.000 & 0.00141881 & 155 & 1.000 & 0.00361169 & 136 & 1.000 \\
  $10^{-6}$ & 0.000392395 & 346 & 1.000 & 0.000394516 & 336 & 1.000 & 0.00227811 & 177 & 1.000 \\
  $10^{-7}$ & 0.00013736 & 795 & 1.000 & 0.000137495 & 798 & 1.000 & 0.00198749 & 388 & 1.000 \\
  $10^{-8}$ & 0.0000702091 & 1705 & 1.000 & 0.0000703217 & 1707 & 1.000 & 0.00190587 & 744 & 1.000 \\
  \hline

  \multicolumn{4}{l||}{$p_\Theta=10^{0}$} &
  \multicolumn{3}{c||}{} &
  \multicolumn{3}{c}{} \\[1mm]
  \hline
  $\gamma$ & gap & iter. & $\lambda$ &
  gap & iter. & $\lambda$ &
  gap & iter. & $\lambda$ \\[1mm]
  \hline
  $10^{-3}$ & 0.199371 & \textcolor{red}{242} & 0.237 & 0.199712 & \textcolor{red}{219} & 0.235 & 0.0299342 & 118 & 1.000 \\
  $10^{-4}$ & 0.191343 & \textcolor{red}{206} & 0.233 & 0.190738 & \textcolor{red}{236} & 0.232 & 0.00935674 & 123 & 1.000 \\
  $10^{-5}$ & 0.178722 & \textcolor{red}{305} & 0.243 & 0.176123 & \textcolor{red}{308} & 0.244 & 0.00415499 & 130 & 1.000 \\
  $10^{-6}$ & 0.175895 & \textcolor{red}{320} & 0.246 & 0.172562 & \textcolor{red}{260} & 0.248 & 0.00301122 & 171 & 1.000 \\
  $10^{-7}$ & 0.175563 & \textcolor{red}{252} & 0.246 & 0.172146 & \textcolor{red}{264} & 0.248 & 0.00279944 & 330 & 1.000 \\
  $10^{-8}$ & 0.175529 & \textcolor{red}{255} & 0.246 & 0.172099 & \textcolor{red}{271} & 0.248 & 0.00275838 & 623 & 1.000 \\
  \hline

  \multicolumn{4}{l||}{$p_\Theta=10^{1}$} &
  \multicolumn{3}{c||}{} &
  \multicolumn{3}{c}{} \\[1mm]
  \hline
  $\gamma$ & gap & iter. & $\lambda$ &
  gap & iter. & $\lambda$ &
  gap & iter. & $\lambda$ \\[1mm]
  \hline
  $10^{-3}$ & 0.215462 & \textcolor{red}{307} & 0.203 & 0.214732 & \textcolor{red}{222} & 0.202 & 0.0301698 & 118 & 1.000 \\
  $10^{-4}$ & 0.149009 & \textcolor{red}{316} & 0.264 & 0.275714 & \textcolor{red}{197} & 0.093 & 0.00946294 & 122 & 1.000 \\
  $10^{-5}$ & 0.127235 & \textcolor{red}{335} & 0.286 & 0.125879 & \textcolor{red}{277} & 0.286 & 0.00459859 & 131 & 1.000 \\
  $10^{-6}$ & 0.123878 & \textcolor{red}{330} & 0.290 & 0.122442 & \textcolor{red}{258} & 0.289 & 0.0037075 & 200 & 1.000 \\
  $10^{-7}$ & 0.123527 & \textcolor{red}{375} & 0.290 & 0.122072 & \textcolor{red}{288} & 0.290 & 0.00338479 & 352 & 1.000 \\
  $10^{-8}$ & 0.123482 & \textcolor{red}{303} & 0.290 & 0.122041 & \textcolor{red}{245} & 0.290 & 0.00338458 & 671 & 1.000 \\
  \hline

  \multicolumn{4}{l||}{$p_\Theta=10^{2}$} &
  \multicolumn{3}{c||}{} &
  \multicolumn{3}{c}{} \\[1mm]
  \hline
  $\gamma$ & gap & iter. & $\lambda$ &
  gap & iter. & $\lambda$ &
  gap & iter. & $\lambda$ \\[1mm]
  \hline
  $10^{-3}$ & 0.286254 & \textcolor{red}{233} & 0.079 & 0.28639 & \textcolor{red}{227} & 0.077 & 0.0313216 & 118 & 1.000 \\
  $10^{-4}$ & 0.253169 & \textcolor{red}{286} & 0.160 & 0.252418 & \textcolor{red}{292} & 0.158 & 0.00926847 & 126 & 1.000 \\
  $10^{-5}$ & 0.286129 & \textcolor{red}{201} & 0.079 & 0.25055 & \textcolor{red}{302} & 0.159 & 0.00498379 & 153 & 1.000 \\
  $10^{-6}$ & 0.286127 & \textcolor{red}{182} & 0.079 & 0.250349 & \textcolor{red}{295} & 0.159 & 0.00434463 & 246 & 1.000 \\
  $10^{-7}$ & 0.286127 & \textcolor{red}{172} & 0.079 & 0.250332 & \textcolor{red}{288} & 0.159 & 0.00414089 & 380 & 1.000 \\
  $10^{-8}$ & 0.286127 & \textcolor{red}{173} & 0.079 & 0.250328 & \textcolor{red}{289} & 0.159 & 0.00418222 & 566 & 1.000 \\
  \hline

  \multicolumn{4}{l||}{$p_\Theta=10^{3}$} &
  \multicolumn{3}{c||}{} &
  \multicolumn{3}{c}{} \\[1mm]
  \hline
  $\gamma$ & gap & iter. & $\lambda$ &
  gap & iter. & $\lambda$ &
  gap & iter. & $\lambda$ \\[1mm]
  \hline
  $10^{-3}$ & 0.287682 & \textcolor{red}{204} & 0.077 & 0.199417 & \textcolor{red}{369} & 0.238 & 0.0614556 & 124 & 1.000 \\
  $10^{-4}$ & 0.277586 & \textcolor{red}{509} & 0.129 & 0.210031 & \textcolor{red}{282} & 0.221 & 0.0261318 & 163 & 1.000 \\
  $10^{-5}$ & 0.263321 & \textcolor{red}{383} & 0.154 & 0.261299 & \textcolor{red}{195} & 0.149 & 0.0161305 & 240 & 1.000 \\
  $10^{-6}$ & 0.287578 & \textcolor{red}{199} & 0.077 & 0.261175 & \textcolor{red}{211} & 0.149 & 0.0155221 & 361 & 1.000 \\
  $10^{-7}$ & 0.287578 & \textcolor{red}{188} & 0.077 & 0.26116 & \textcolor{red}{182} & 0.149 & 0.059845 & \textcolor{red}{567} & 0.941 \\
  $10^{-8}$ & 0.287578 & \textcolor{red}{194} & 0.077 & 0.261158 & \textcolor{red}{197} & 0.149 & 0.059578 & \textcolor{red}{902} & 0.942 \\
  \hline\hline
\end{tabular}
\label{tab:sensitivity_study_ptheta_cont}
\end{Table}

\newpage

\begin{Table}[H]
\centering
\scriptsize
\setlength{\tabcolsep}{3pt}
\caption{
Residual gap and total Newton iteration count for the self-contact-within-a-box benchmark.
Results are listed for different combinations of the third-medium scaling parameter $\gamma$ and the regularization parameter $\alpha_r$.
The penalty parameter is fixed to $p_\Theta=10^{-1}$.
}
\begin{tabular}{c|c|c||c|c|c||c|c|c||c|c|c}
  \hline\hline
  \multicolumn{3}{c||}{T$_2^u$T$_2^\Theta$} &
  \multicolumn{3}{c||}{T$_2^u$T$_1^\Theta$} &
  \multicolumn{3}{c||}{T$_1^u$T$_1^\Theta$} &
  \multicolumn{3}{c}{T$_1^u$T$_0^{\Theta,d}$} \\[1mm]
  \hline

  \multicolumn{3}{l||}{$\alpha_r=10^{-2}$} &
  \multicolumn{3}{c||}{} &
  \multicolumn{3}{c||}{} &
  \multicolumn{3}{c}{} \\[1mm]
  \hline
  $\gamma$ & gap & iter. &
  $\gamma$ & gap & iter. &
  $\gamma$ & gap & iter. &
  $\gamma$ & gap & iter. \\[1mm]
  \hline
  $10^{-3}$ & 0.031649 & 120 & $10^{-3}$ & 0.0315106 & 120 & $10^{-3}$ & 0.0312422 & 120 & $10^{-3}$ & 0.0297947 & 118 \\
  $10^{-4}$ & 0.0084414 & 120 & $10^{-4}$ & 0.00838241 & 120 & $10^{-4}$ & 0.00943054 & 123 & $10^{-4}$ & 0.00910682 & 125 \\
  $10^{-5}$ & 0.00194973 & 134 & $10^{-5}$ & 0.00194097 & 133 & $10^{-5}$ & 0.00402415 & 139 & $10^{-5}$ & 0.00351142 & 137 \\
  $10^{-6}$ & 0.000515392 & 228 & $10^{-6}$ & 0.000516935 & 234 & $10^{-6}$ & 0.00299823 & 270 & $10^{-6}$ & 0.00214174 & 194 \\
  $10^{-7}$ & 0.000176805 & 538 & $10^{-7}$ & 0.00017887 & 571 & $10^{-7}$ & 0.00276324 & 535 & $10^{-7}$ & 0.00180956 & 405 \\
  $10^{-8}$ & 0.0000862862 & 1167 & $10^{-8}$ & 0.0000888917 & 1281 & $10^{-8}$ & 0.00269713 & 993 & $10^{-8}$ & 0.00173047 & 891 \\
  \hline

  \multicolumn{3}{l||}{$\alpha_r=10^{-1}$} &
  \multicolumn{3}{c||}{} &
  \multicolumn{3}{c||}{} &
  \multicolumn{3}{c}{} \\[1mm]
  \hline
  $\gamma$ & gap & iter. &
  $\gamma$ & gap & iter. &
  $\gamma$ & gap & iter. &
  $\gamma$ & gap & iter. \\[1mm]
  \hline
  $10^{-3}$ & 0.0264975 & 119 & $10^{-3}$ & 0.0264158 & 119 & $10^{-3}$ & 0.0266828 & 119 & $10^{-3}$ & 0.0297947 & 118 \\
  $10^{-4}$ & 0.00643342 & 128 & $10^{-4}$ & 0.00640789 & 127 & $10^{-4}$ & 0.00800657 & 128 & $10^{-4}$ & 0.00910682 & 125 \\
  $10^{-5}$ & 0.00147664 & 148 & $10^{-5}$ & 0.0014734 & 148 & $10^{-5}$ & 0.00383879 & 158 & $10^{-5}$ & 0.00351142 & 137 \\
  $10^{-6}$ & 0.000419329 & 310 & $10^{-6}$ & 0.000415996 & 319 & $10^{-6}$ & 0.00296105 & 346 & $10^{-6}$ & 0.00214174 & 194 \\
  $10^{-7}$ & 0.000144945 & 742 & $10^{-7}$ & 0.000145813 & 745 & $10^{-7}$ & 0.00275866 & 690 & $10^{-7}$ & 0.00180956 & 405 \\
  $10^{-8}$ & 0.0000733487 & 1628 & $10^{-8}$ & 0.0000742735 & 1612 & $10^{-8}$ & 0.00270738 & 1253 & $10^{-8}$ & 0.00173047 & 891 \\
  \hline

  \multicolumn{3}{l||}{$\alpha_r=10^{0}$} &
  \multicolumn{3}{c||}{} &
  \multicolumn{3}{c||}{} &
  \multicolumn{3}{c}{} \\[1mm]
  \hline
  $\gamma$ & gap & iter. &
  $\gamma$ & gap & iter. &
  $\gamma$ & gap & iter. &
  $\gamma$ & gap & iter. \\[1mm]
  \hline
  $10^{-3}$ & 0.0251904 & 119 & $10^{-3}$ & 0.0251813 & 119 & $10^{-3}$ & 0.0255675 & 119 & $10^{-3}$ & 0.0297947 & 118 \\
  $10^{-4}$ & 0.00600028 & 128 & $10^{-4}$ & 0.00599759 & 128 & $10^{-4}$ & 0.0077305 & 129 & $10^{-4}$ & 0.00910682 & 125 \\
  $10^{-5}$ & 0.00141917 & 155 & $10^{-5}$ & 0.00141881 & 155 & $10^{-5}$ & 0.00382607 & 171 & $10^{-5}$ & 0.00351142 & 137 \\
  $10^{-6}$ & 0.000392395 & 346 & $10^{-6}$ & 0.000394516 & 336 & $10^{-6}$ & 0.00299808 & 366 & $10^{-6}$ & 0.00214174 & 194 \\
  $10^{-7}$ & 0.00013736 & 795 & $10^{-7}$ & 0.000137495 & 798 & $10^{-7}$ & 0.00280456 & 725 & $10^{-7}$ & 0.00180956 & 405 \\
  $10^{-8}$ & 0.0000702091 & 1705 & $10^{-8}$ & 0.0000703217 & 1707 & $10^{-8}$ & 0.00274656 & 1308 & $10^{-8}$ & 0.00173047 & 891 \\
  \hline

  \multicolumn{3}{l||}{$\alpha_r=10^{1}$} &
  \multicolumn{3}{c||}{} &
  \multicolumn{3}{c||}{} &
  \multicolumn{3}{c}{} \\[1mm]
  \hline
  $\gamma$ & gap & iter. &
  $\gamma$ & gap & iter. &
  $\gamma$ & gap & iter. &
  $\gamma$ & gap & iter. \\[1mm]
  \hline
  $10^{-3}$ & 0.0250555 & 118 & $10^{-3}$ & 0.0250546 & 118 & $10^{-3}$ & 0.0254551 & 118 & $10^{-3}$ & 0.0297947 & 118 \\
  $10^{-4}$ & 0.00595654 & 128 & $10^{-4}$ & 0.00595627 & 127 & $10^{-4}$ & 0.00770511 & 128 & $10^{-4}$ & 0.00910682 & 125 \\
  $10^{-5}$ & 0.00140966 & 154 & $10^{-5}$ & 0.00140962 & 154 & $10^{-5}$ & 0.00384868 & 191 & $10^{-5}$ & 0.00351142 & 137 \\
  $10^{-6}$ & 0.000393792 & 339 & $10^{-6}$ & 0.000394525 & 338 & $10^{-6}$ & 0.00299173 & 366 & $10^{-6}$ & 0.00214174 & 194 \\
  $10^{-7}$ & 0.00013654 & 807 & $10^{-7}$ & 0.000136633 & 795 & $10^{-7}$ & 0.00281297 & 741 & $10^{-7}$ & 0.00180956 & 405 \\
  $10^{-8}$ & 0.0000698996 & 1735 & $10^{-8}$ & 0.0000699086 & 1711 & $10^{-8}$ & 0.00275783 & 1347 & $10^{-8}$ & 0.00173047 & 891 \\
  \hline

  \multicolumn{3}{l||}{$\alpha_r=10^{2}$} &
  \multicolumn{3}{c||}{} &
  \multicolumn{3}{c||}{} &
  \multicolumn{3}{c}{} \\[1mm]
  \hline
  $\gamma$ & gap & iter. &
  $\gamma$ & gap & iter. &
  $\gamma$ & gap & iter. &
  $\gamma$ & gap & iter. \\[1mm]
  \hline
  $10^{-3}$ & 0.0250421 & 118 & $10^{-3}$ & 0.025042 & 118 & $10^{-3}$ & 0.0254439 & 118 & $10^{-3}$ & 0.0297947 & 118 \\
  $10^{-4}$ & 0.00595223 & 128 & $10^{-4}$ & 0.0059522 & 127 & $10^{-4}$ & 0.00770265 & 128 & $10^{-4}$ & 0.00910682 & 125 \\
  $10^{-5}$ & 0.00140872 & 154 & $10^{-5}$ & 0.00140871 & 154 & $10^{-5}$ & 0.00381932 & 177 & $10^{-5}$ & 0.00351142 & 137 \\
  $10^{-6}$ & 0.000393792 & 340 & $10^{-6}$ & 0.000393943 & 339 & $10^{-6}$ & 0.00299839 & 368 & $10^{-6}$ & 0.00214174 & 194 \\
  $10^{-7}$ & 0.000136544 & 797 & $10^{-7}$ & 0.000136505 & 802 & $10^{-7}$ & 0.0028137 & 737 & $10^{-7}$ & 0.00180956 & 405 \\
  $10^{-8}$ & 0.0000698444 & 1751 & $10^{-8}$ & 0.0000698244 & 1707 & $10^{-8}$ & 0.00276193 & 1341 & $10^{-8}$ & 0.00173047 & 891 \\
  \hline

  \multicolumn{3}{l||}{$\alpha_r=10^{3}$} &
  \multicolumn{3}{c||}{} &
  \multicolumn{3}{c||}{} &
  \multicolumn{3}{c}{} \\[1mm]
  \hline
  $\gamma$ & gap & iter. &
  $\gamma$ & gap & iter. &
  $\gamma$ & gap & iter. &
  $\gamma$ & gap & iter. \\[1mm]
  \hline
  $10^{-3}$ & 0.0250408 & 118 & $10^{-3}$ & 0.0250408 & 118 & $10^{-3}$ & 0.0254428 & 118 & $10^{-3}$ & 0.0297947 & 118 \\
  $10^{-4}$ & 0.0059518 & 128 & $10^{-4}$ & 0.0059518 & 127 & $10^{-4}$ & 0.00770241 & 128 & $10^{-4}$ & 0.00910682 & 125 \\
  $10^{-5}$ & 0.00140862 & 154 & $10^{-5}$ & 0.00140862 & 154 & $10^{-5}$ & 0.00384226 & 174 & $10^{-5}$ & 0.00351142 & 137 \\
  $10^{-6}$ & 0.000393987 & 341 & $10^{-6}$ & 0.000394244 & 338 & $10^{-6}$ & 0.00299542 & 366 & $10^{-6}$ & 0.00214174 & 194 \\
  $10^{-7}$ & 0.000136242 & 817 & $10^{-7}$ & 0.000136496 & 805 & $10^{-7}$ & 0.00281715 & 727 & $10^{-7}$ & 0.00180956 & 405 \\
  $10^{-8}$ & 0.0000698494 & 1711 & $10^{-8}$ & 0.0000698357 & 1708 & $10^{-8}$ & 0.00276403 & 1344 & $10^{-8}$ & 0.00173047 & 891 \\
  \hline

  \multicolumn{3}{l||}{$\alpha_r=10^{4}$} &
  \multicolumn{3}{c||}{} &
  \multicolumn{3}{c||}{} &
  \multicolumn{3}{c}{} \\[1mm]
  \hline
  $\gamma$ & gap & iter. &
  $\gamma$ & gap & iter. &
  $\gamma$ & gap & iter. &
  $\gamma$ & gap & iter. \\[1mm]
  \hline
  $10^{-3}$ & 0.0250406 & 119 & $10^{-3}$ & 0.0250406 & 118 & $10^{-3}$ & 0.0254427 & 118 & $10^{-3}$ & 0.0297947 & 118 \\
  $10^{-4}$ & 0.00595176 & 129 & $10^{-4}$ & 0.00595176 & 127 & $10^{-4}$ & 0.00770239 & 128 & $10^{-4}$ & 0.00910682 & 125 \\
  $10^{-5}$ & 0.00140861 & 155 & $10^{-5}$ & 0.00140861 & 154 & $10^{-5}$ & 0.00384226 & 174 & $10^{-5}$ & 0.00351142 & 137 \\
  $10^{-6}$ & 0.000393984 & 342 & $10^{-6}$ & 0.000394242 & 338 & $10^{-6}$ & 0.00299542 & 366 & $10^{-6}$ & 0.00214174 & 194 \\
  $10^{-7}$ & 0.000136463 & 819 & $10^{-7}$ & 0.000136498 & 805 & $10^{-7}$ & 0.00281716 & 727 & $10^{-7}$ & 0.00180956 & 405 \\
  $10^{-8}$ & 0.0000698491 & 1712 & $10^{-8}$ & 0.0000698371 & 1707 & $10^{-8}$ & 0.00276403 & 1344 & $10^{-8}$ & 0.00173047 & 891 \\
  \hline\hline
\end{tabular}
\label{tab:sensitivity_study_alpha_r_cont}
\end{Table}

\begin{Table}[H]
\centering
\scriptsize
\setlength{\tabcolsep}{3pt}
\caption{
Residual gap, total Newton iteration count and achieved load parameter for the self-contact-within-a-box benchmark using discontinuous auxiliary-field interpolations.
Results are listed for different combinations of the third-medium scaling parameter $\gamma$ and the penalty parameter $p_\Theta$, while $\alpha_r=1$.
Red iteration counts denote diverged parameter sets with load parameter $\lambda<0.95$, whereas orange iteration counts denote simulations that diverged at $\lambda\geq0.95$.
}
\begin{tabular}{c|c|c|c||c|c|c||c|c|c||c|c|c}
  \hline\hline
  \multicolumn{4}{c||}{T$_2^u$T$_1^{\Theta,d}$} &
  \multicolumn{3}{c||}{T$_2^u$T$_0^{\Theta,d}$} &
  \multicolumn{3}{c||}{T$_1^u$T$_1^{\Theta,d}$} &
  \multicolumn{3}{c}{T$_1^u$T$_0^{\Theta,d}$} \\[1mm]
  \hline

  \multicolumn{4}{l||}{$p_\Theta=10^{-2}$} &
  \multicolumn{3}{c||}{} &
  \multicolumn{3}{c||}{} &
  \multicolumn{3}{c}{} \\[1mm]
  \hline
  $\gamma$ & gap & iter. & $\lambda$ &
  gap & iter. & $\lambda$ &
  gap & iter. & $\lambda$ &
  gap & iter. & $\lambda$ \\[1mm]
  \hline
  $10^{-3}$ & 0.0298569 & 118 & 1.000 & 0.0298569 & 118 & 1.000 & 0.0297085 & 118 & 1.000 & 0.0297246 & 118 & 1.000 \\
  $10^{-4}$ & 0.00802477 & 126 & 1.000 & 0.00802477 & 126 & 1.000 & 0.00824478 & 174 & 1.000 & 0.00876471 & 126 & 1.000 \\
  $10^{-5}$ & 0.00196045 & 147 & 1.000 & 0.00196045 & 147 & 1.000 & 0.00265286 & \textcolor{red}{365} & 0.710 & 0.00242643 & 160 & 1.000 \\
  $10^{-6}$ & 0.000522379 & \textcolor{orange}{476} & 0.983 & 0.000522379 & \textcolor{orange}{476} & 0.983 & 0.000992186 & \textcolor{red}{315} & 0.487 & 0.000860014 & \textcolor{red}{439} & 0.891 \\
  $10^{-7}$ & 0.000146887 & \textcolor{red}{376} & 0.613 & 0.000146887 & \textcolor{red}{376} & 0.613 & 0.000336383 & \textcolor{red}{325} & 0.426 & 0.00055824 & \textcolor{red}{265} & 0.430 \\
  $10^{-8}$ & 0.0000445499 & \textcolor{red}{445} & 0.461 & 0.0000445499 & \textcolor{red}{445} & 0.461 & 0.0000386906 & \textcolor{red}{437} & 0.453 & 0.000331303 & \textcolor{red}{278} & 0.413 \\
  \hline

  \multicolumn{4}{l||}{$p_\Theta=10^{-1}$} &
  \multicolumn{3}{c||}{} &
  \multicolumn{3}{c||}{} &
  \multicolumn{3}{c}{} \\[1mm]
  \hline
  $\gamma$ & gap & iter. & $\lambda$ &
  gap & iter. & $\lambda$ &
  gap & iter. & $\lambda$ &
  gap & iter. & $\lambda$ \\[1mm]
  \hline
  $10^{-3}$ & 0.0298864 & 118 & 1.000 & 0.0298864 & 118 & 1.000 & 0.0297085 & 118 & 1.000 & 0.0297947 & 118 & 1.000 \\
  $10^{-4}$ & 0.00805854 & 126 & 1.000 & 0.00805854 & 126 & 1.000 & 0.00824478 & 174 & 1.000 & 0.00910682 & 125 & 1.000 \\
  $10^{-5}$ & 0.00199623 & 148 & 1.000 & 0.00199623 & 148 & 1.000 & 0.00265286 & \textcolor{red}{365} & 0.710 & 0.00351142 & 137 & 1.000 \\
  $10^{-6}$ & 0.000545862 & \textcolor{red}{276} & 0.690 & 0.000545862 & \textcolor{red}{276} & 0.690 & 0.000992186 & \textcolor{red}{315} & 0.487 & 0.00214174 & 194 & 1.000 \\
  $10^{-7}$ & 0.000134083 & \textcolor{red}{374} & 0.581 & 0.000134083 & \textcolor{red}{374} & 0.581 & 0.000336383 & \textcolor{red}{325} & 0.426 & 0.00180956 & 405 & 1.000 \\
  $10^{-8}$ & 0.000240972 & \textcolor{red}{276} & 0.353 & 0.000240972 & \textcolor{red}{276} & 0.353 & 0.0000386906 & \textcolor{red}{437} & 0.453 & 0.00173047 & 891 & 1.000 \\
  \hline

  \multicolumn{4}{l||}{$p_\Theta=10^{0}$} &
  \multicolumn{3}{c||}{} &
  \multicolumn{3}{c||}{} &
  \multicolumn{3}{c}{} \\[1mm]
  \hline
  $\gamma$ & gap & iter. & $\lambda$ &
  gap & iter. & $\lambda$ &
  gap & iter. & $\lambda$ &
  gap & iter. & $\lambda$ \\[1mm]
  \hline
  $10^{-3}$ & 0.0300636 & 118 & 1.000 & 0.0300636 & 118 & 1.000 & 0.0297085 & 118 & 1.000 & 0.0301713 & 118 & 1.000 \\
  $10^{-4}$ & 0.00833051 & 127 & 1.000 & 0.00833051 & 127 & 1.000 & 0.00824478 & 174 & 1.000 & 0.00951916 & 123 & 1.000 \\
  $10^{-5}$ & 0.00211543 & 229 & 1.000 & 0.00211543 & 229 & 1.000 & 0.00265286 & \textcolor{red}{365} & 0.710 & 0.00417045 & 131 & 1.000 \\
  $10^{-6}$ & 0.00121318 & \textcolor{red}{347} & 0.438 & 0.00121318 & \textcolor{red}{347} & 0.438 & 0.000992186 & \textcolor{red}{315} & 0.487 & 0.0030107 & 194 & 1.000 \\
  $10^{-7}$ & 0.000211673 & \textcolor{red}{375} & 0.438 & 0.000211673 & \textcolor{red}{375} & 0.438 & 0.000336383 & \textcolor{red}{325} & 0.426 & 0.00275762 & 377 & 1.000 \\
  $10^{-8}$ & 0.0000677622 & \textcolor{red}{289} & 0.385 & 0.0000677622 & \textcolor{red}{289} & 0.385 & 0.0000386906 & \textcolor{red}{437} & 0.453 & 0.00268435 & 752 & 1.000 \\
  \hline

  \multicolumn{4}{l||}{$p_\Theta=10^{1}$} &
  \multicolumn{3}{c||}{} &
  \multicolumn{3}{c||}{} &
  \multicolumn{3}{c}{} \\[1mm]
  \hline
  $\gamma$ & gap & iter. & $\lambda$ &
  gap & iter. & $\lambda$ &
  gap & iter. & $\lambda$ &
  gap & iter. & $\lambda$ \\[1mm]
  \hline
  $10^{-3}$ & 0.0303078 & 118 & 1.000 & 0.0303078 & 118 & 1.000 & 0.0297085 & 118 & 1.000 & 0.0316236 & 118 & 1.000 \\
  $10^{-4}$ & 0.00857492 & 128 & 1.000 & 0.00857492 & 128 & 1.000 & 0.00824478 & 174 & 1.000 & 0.0113235 & 124 & 1.000 \\
  $10^{-5}$ & 0.0021616 & 279 & 1.000 & 0.0021616 & 279 & 1.000 & 0.00265286 & \textcolor{red}{365} & 0.710 & 0.00639403 & 146 & 1.000 \\
  $10^{-6}$ & 0.000734647 & \textcolor{red}{343} & 0.547 & 0.000734647 & \textcolor{red}{343} & 0.547 & 0.000992186 & \textcolor{red}{315} & 0.487 & 0.00527773 & 239 & 1.000 \\
  $10^{-7}$ & 0.000282756 & \textcolor{red}{405} & 0.487 & 0.000282756 & \textcolor{red}{405} & 0.487 & 0.000336383 & \textcolor{red}{325} & 0.426 & 0.00504656 & 480 & 1.000 \\
  $10^{-8}$ & 0.000499722 & \textcolor{red}{257} & 0.349 & 0.000499722 & \textcolor{red}{257} & 0.349 & 0.0000386906 & \textcolor{red}{437} & 0.453 & 0.00498715 & 791 & 1.000 \\
  \hline

  \multicolumn{4}{l||}{$p_\Theta=10^{2}$} &
  \multicolumn{3}{c||}{} &
  \multicolumn{3}{c||}{} &
  \multicolumn{3}{c}{} \\[1mm]
  \hline
  $\gamma$ & gap & iter. & $\lambda$ &
  gap & iter. & $\lambda$ &
  gap & iter. & $\lambda$ &
  gap & iter. & $\lambda$ \\[1mm]
  \hline
  $10^{-3}$ & 0.0303545 & 118 & 1.000 & 0.0303545 & 118 & 1.000 & 0.0297085 & 118 & 1.000 & 0.0363537 & 119 & 1.000 \\
  $10^{-4}$ & 0.00861313 & 135 & 1.000 & 0.00861313 & 135 & 1.000 & 0.00824478 & 174 & 1.000 & 0.0181991 & 130 & 1.000 \\
  $10^{-5}$ & 0.00216338 & \textcolor{orange}{434} & 0.960 & 0.00216338 & \textcolor{orange}{434} & 0.960 & 0.00265286 & \textcolor{red}{365} & 0.710 & 0.0139122 & 180 & 1.000 \\
  $10^{-6}$ & 0.000742145 & \textcolor{red}{369} & 0.536 & 0.000742145 & \textcolor{red}{369} & 0.536 & 0.000992186 & \textcolor{red}{315} & 0.487 & 0.0129623 & 296 & 1.000 \\
  $10^{-7}$ & 0.000187948 & \textcolor{red}{336} & 0.448 & 0.000187948 & \textcolor{red}{336} & 0.448 & 0.000336383 & \textcolor{red}{325} & 0.426 & 0.012751 & 478 & 1.000 \\
  $10^{-8}$ & 0.000099741 & \textcolor{red}{385} & 0.440 & 0.000099741 & \textcolor{red}{385} & 0.440 & 0.0000386906 & \textcolor{red}{437} & 0.453 & 0.0127011 & 682 & 1.000 \\
  \hline

  \multicolumn{4}{l||}{$p_\Theta=10^{3}$} &
  \multicolumn{3}{c||}{} &
  \multicolumn{3}{c||}{} &
  \multicolumn{3}{c}{} \\[1mm]
  \hline
  $\gamma$ & gap & iter. & $\lambda$ &
  gap & iter. & $\lambda$ &
  gap & iter. & $\lambda$ &
  gap & iter. & $\lambda$ \\[1mm]
  \hline
  $10^{-3}$ & 0.0303416 & 118 & 1.000 & 0.0303416 & 118 & 1.000 & 0.0297085 & 118 & 1.000 & 0.0552937 & 119 & 1.000 \\
  $10^{-4}$ & 0.00860784 & 134 & 1.000 & 0.00860784 & 134 & 1.000 & 0.00824478 & 174 & 1.000 & 0.0379871 & 170 & 1.000 \\
  $10^{-5}$ & 0.00231965 & \textcolor{red}{389} & 0.872 & 0.00231965 & \textcolor{red}{389} & 0.872 & 0.00265286 & \textcolor{red}{365} & 0.710 & 0.034976 & 229 & 1.000 \\
  $10^{-6}$ & 0.0011114 & \textcolor{red}{243} & 0.471 & 0.0011114 & \textcolor{red}{243} & 0.471 & 0.000992186 & \textcolor{red}{315} & 0.487 & 0.0340599 & 296 & 1.000 \\
  $10^{-7}$ & 0.000190877 & \textcolor{red}{286} & 0.456 & 0.000190877 & \textcolor{red}{286} & 0.456 & 0.000336383 & \textcolor{red}{330} & 0.426 & 0.033939 & 381 & 1.000 \\
  $10^{-8}$ & 0.0000382093 & \textcolor{red}{458} & 0.440 & 0.0000382093 & \textcolor{red}{458} & 0.440 & 0.0000386906 & \textcolor{red}{441} & 0.453 & 0.0336688 & 548 & 1.000 \\
  \hline

  \multicolumn{4}{l||}{$p_\Theta=10^{4}$} &
  \multicolumn{3}{c||}{} &
  \multicolumn{3}{c||}{} &
  \multicolumn{3}{c}{} \\[1mm]
  \hline
  $\gamma$ & gap & iter. & $\lambda$ &
  gap & iter. & $\lambda$ &
  gap & iter. & $\lambda$ &
  gap & iter. & $\lambda$ \\[1mm]
  \hline
  $10^{-3}$ & 0.0303393 & 118 & 1.000 & 0.0303393 & 118 & 1.000 & 0.0297085 & 118 & 1.000 & 0.0809935 & 188 & 1.000 \\
  $10^{-4}$ & 0.00860684 & 133 & 1.000 & 0.00860684 & 133 & 1.000 & 0.00824478 & 174 & 1.000 & 0.13106 & \textcolor{red}{345} & 0.724 \\
  $10^{-5}$ & 0.00232095 & \textcolor{red}{452} & 0.871 & 0.00232095 & \textcolor{red}{452} & 0.871 & 0.00265286 & \textcolor{red}{364} & 0.710 & 0.12846 & \textcolor{red}{441} & 0.720 \\
  $10^{-6}$ & 0.000826532 & \textcolor{red}{411} & 0.674 & 0.000826532 & \textcolor{red}{411} & 0.674 & 0.000992186 & \textcolor{red}{315} & 0.487 & 0.128158 & \textcolor{red}{584} & 0.719 \\
  $10^{-7}$ & 0.000218183 & \textcolor{red}{283} & 0.452 & 0.000218183 & \textcolor{red}{283} & 0.452 & 0.000336383 & \textcolor{red}{330} & 0.426 & 0.128099 & \textcolor{red}{776} & 0.719 \\
  $10^{-8}$ & 0.00032454 & \textcolor{red}{240} & 0.351 & 0.00032454 & \textcolor{red}{240} & 0.351 & 0.0000386906 & \textcolor{red}{441} & 0.453 & 0.128091 & \textcolor{red}{1253} & 0.719 \\
  \hline\hline
\end{tabular}
\label{tab:sensitivity_study_ptheta_disc}
\end{Table}

\end{appendix}

\end{document}